\documentclass[12pt,reqno]{amsart}
\usepackage{xifthen}
\usepackage{bm}
\usepackage{subfig}
\usepackage[normalem]{ulem}
\usepackage{pkg}
\usepackage{epsfig}
\usepackage{eucal}

\numberwithin{equation}{section}

\newcommand{\AthreyaNey}{MR2047480}

\newcommand{\Bill}{MR1324786}
\newcommand{\BinghamDoney}{MR0362525}
\newcommand{\BolthausenSznitman}{MR2023130}
\newcommand{\BNT}{MR1177900}
\newcommand{\BVCvx}{MR2061575}
\newcommand{\CroydonFK}{arXiv:1203.4078}

\newcommand{\DZ}{MR1619036}

\newcommand{\EK}{MR838085}
\newcommand{\Faraud}{MR2754802}

\newcommand{\FellerII}{MR0270403}
\newcommand{\GWLeaves}{AOP0810.012R1A0}
\newcommand{\Harris}{MR1991122}
\newcommand{\HJ}{MR1084815}
\newcommand{\HuShi}{MR2299718}
\newcommand{\Kallenberg}{MR1876169}
\newcommand{\KLPP}{MR1601737}
\newcommand{\KSBRW}{MR1849369}
\newcommand{\KS}{MR1121940}
\newcommand{\KSK}{MR0407981}
\newcommand{\KSMGW}{MR0198552}
\newcommand{\Liu}{MR1741808}
\newcommand{\LPRWRE}{MR1143414}

\newcommand{\LPPErgodic}{MR1336708}
\newcommand{\LPPBiased}{MR1410689}
\newcommand{\Lyons}{MR1062053}

\newcommand{\NeyVMoments}{MR1970272}
\newcommand{\Petersen}{MR1073173}
\newcommand{\Petrov}{MR0388499}
\newcommand{\PPPolar}{MR1349163}
\newcommand{\PZ}{MR2365486}

\newcommand{\Williams}{MR1155402}
\newcommand{\Zeitouni}{MR2071631}

\newcommand{\hypirred}{\text{\textup{(H1)}}}
\newcommand{\hypks}{\text{\textup{(H2)}}}
\newcommand{\hypmmt}[1]{\text{\textup{(H3$^{#1}$)}}}
\newcommand{\mmt}{p}

\newcommand{\bpi}{{\bm{\pi}}}
\newcommand{\bthet}{\bm{\theta}}
\newcommand{\g}{{\mathbf{g}}}
\newcommand{\infq}{\wh\q}
\newcommand{\q}{\mathbf{q}}

\newcommand{\leftevec}{g}
\newcommand{\lv}{\vec\leftevec}
\newcommand{\lvec}[1]{\leftevec_{#1}}

\newcommand{\rightevec}{e}

\newcommand{\emin}{{\rightevec_{\min}}}

\newcommand{\rv}{{\vec\rightevec}}
\newcommand{\rvec}[1]{{\rightevec_{#1}}}

\newcommand{\cond}[2]{\cC_{#1,#2}}
\newcommand{\desc}{\partial^+}
\newcommand{\leaves}{\cL}
\newcommand{\subtree}[2]{#1^{(#2)}}
\newcommand{\ray}{\xi}
\newcommand{\lne}{\bar\xi}
\newcommand{\rt}{o}
\newcommand{\shift}{\mathfrak{S}}
\newcommand{\trees}{\Om}
\newcommand{\treesray}{{\trees_\downarrow}}
\newcommand{\treesline}{{\trees_\updownarrow}}

\newcommand{\treesrayinf}{\Om^\infty_\downarrow}
\newcommand{\type}{\chi}
\DeclareMathOperator{\br}{br}

\newcommand{\acr}[1]{{\text{\textsf{\textup{#1}}}}}

\newcommand{\infl}[1]{{\acr{I}#1}}
\newcommand{\revers}[1]{{#1\acr{R}}}
\newcommand{\gwt}{\acr{GW}}

\newcommand{\igw}{\infl\gwt}
\newcommand{\igwr}{{\revers\igw}}
\newcommand{\mgw}{{\acr{MGW}}}
\newcommand{\mgwtype}[1]{{\mgw^{#1}}}

\newcommand{\imgw}{{\infl\mgw}}
\newcommand{\imgwr}{{\revers\imgw}}
\newcommand{\Qn}{\acr{Q}}
\newcommand{\Qntype}[1]{\Qn^{#1}} 
\newcommand{\barimgwztype}[1]{\underleftarrow\imgw^{#1}_0}
\newcommand{\infmgw}{{\widehat\mgw}}
\newcommand{\infmgwtype}[1]{\infmgw^{#1}}
\newcommand{\infmgwstar}{{\infmgw_\star}}
\newcommand{\infmgwstartype}[1]{\infmgwstar^{#1}}
\newcommand{\imgwline}{\underline\imgw}

\newcommand{\Znorm}{\mathfrak{Z}}

\newcommand{\cQre}{\bar\cQ}
\newcommand{\rvre}{\bar\rv}
\newcommand{\lvre}{\bar\lv}
\newcommand{\rvecre}[1]{{\bar\rightevec}_{#1}}
\newcommand{\lvecre}[1]{{\bar\leftevec}_{#1}}
\newcommand{\mgwre}{\ol\mgw}
\newcommand{\mgwretype}[1]{\mgwre^{#1}}
\newcommand{\imgwre}{\ol\imgw}
\newcommand{\imgwrre}{\ol\imgwr}
\newcommand{\barimgwztypere}[1]{\underleftarrow{\overline\imgw}^{#1}_0}
\newcommand{\Qnre}{\ol\Qn}

\newcommand{\cts}{\text{\textup{cts}}}
\newcommand{\rw}{\acr{RW}}
\newcommand{\tree}{\cT}
\newcommand{\trw}{\acr{TRW}}
\newcommand{\rwre}{\acr{RWRE}}
\newcommand{\trwre}{\acr{TRWRE}}

\newcommand{\Pinj}{\cP^\mathrm{inj}}
\newcommand{\aT}{\mathrm{\textsf{T}}}

\newcommand{\bs}{\mathbf{s}}
\newcommand{\bt}{\mathbf{t}}
\newcommand{\bb}{\mathbf{b}}
\newcommand{\Thit}{\tau^{\mathrm{hit}}}
\newcommand{\Tret}{\tau^\circlearrowleft}
\newcommand{\nx}{\mathbb{X}}
\newcommand{\nxtime}{\tau^\nx}

\newcommand{\rtime}[2]{{\tau^{#1}(#2)}}
\newcommand{\Troot}[2]{N_{#2}(#1)}
\newcommand{\Trootcoup}[2]{N^\coup_{#2}(#1)}
\newcommand{\visits}{L}
\newcommand{\locT}[2]{L_{#2}(#1)}
\newcommand{\exc}{\mathrm{exc}}
\newcommand{\cent}{\mathrm{cent}}
\newcommand{\Texc}{\tau^{\mathrm{exc}}}
\newcommand{\perc}{\mathrm{perc}}

\newcommand{\coup}{\diamond}
\newcommand{\drift}{\mathbb{D}}
\newcommand{\reflec}{\mathbb{H}}
\newcommand{\inter}{\mathrm{int}}
\newcommand{\badA}{\mathbf{A}}
\newcommand{\badB}{\mathbf{B}}
\newcommand{\badC}{\mathbf{C}}
\newcommand{\badE}{\mathbf{E}}
\newcommand{\badY}{\bm{\Upsilon}}
\newcommand{\badThet}{\mathbf{G}}

\newcommand{\domain}{\mathcal{D}}

\newcommand{\rem}{\mathrm{\textsf{R}}}
\newcommand{\poly}{\mathrm{\textsf{P}}}
\newcommand{\expm}{\mathrm{\textsf{e}}}

\begin{document}

\title[CLT for biased random walk on multi-type Galton--Watson trees]{Central limit theorem for biased random walk \\ on multi-type Galton--Watson trees}

\author[A.\ Dembo]{$^*$Amir Dembo}
\address{$^*$Department of Mathematics, Stanford University
\newline\indent Building 380, Sloan Hall, Stanford, California 94305}
\author[N.\ Sun]{$^\dagger$Nike Sun}
\address{$^{*\dagger}$Department of Statistics, Stanford University
\newline\indent Sequoia Hall, 390 Serra Mall, Stanford, California 94305}

\date{\today}

\subjclass[2010]{Primary 60F05, 60K37; Secondary 60J80, 60G50}
%60F05  	Central limit and other weak theorems
%60K37  	Processes in random environments
%60J80  	Branching processes (Galton--Watson,
%			birth-and-death, etc.)
%60G50  	Sums of independent random variables;
%			random walks

\keywords{Multi-type Galton--Watson tree, biased random walk, central limit theorem, random walk with random environment.}

\thanks{$^{*\dagger}$Research partially supported by NSF grants DMS-0806211 and DMS-1106627. $^\dagger$Research partially supported by NSF grant DMS-0502385 and Department of Defense NDSEG Fellowship.}

\begin{abstract}
Let $\tree$ be a rooted supercritical multi-type Galton--Watson (MGW) tree with types coming from a finite alphabet, conditioned to non-extinction. The $\lm$-biased random walk $(X_t)_{t\ge0}$ on $\tree$ is the nearest-neighbor random walk which, when at a vertex $v$ with $d_v$ offspring, moves closer to the root with probability $\lm/(\lm+d_v)$, and to each of the offspring with probability $1/(\lm+d_v)$. This walk is recurrent for $\lm\ge\rho$ and transient for $0\le\lm<\rho$, with $\rho$ the Perron--Frobenius eigenvalue for the (assumed) irreducible matrix of expected offspring numbers. Subject to finite moments of order $\mmt>4$ for the offspring distributions, we prove the following quenched CLT for $\lm$-biased random walk at the critical value $\lm=\rho$: for almost every $\tree$, the process $\abs{X_{\flr{nt}}}/\sqrt{n}$ converges in law as $n\to\infty$ to a reflected Brownian motion rescaled by an explicit constant. This result was proved under some stronger assumptions by Peres--Zeitouni (2008) for single-type Galton--Watson trees. Following their approach, our proof is based on a new explicit description of a reversing measure for the walk from the point of view of the particle (generalizing the measure constructed in the single-type setting by Peres--Zeitouni), and the construction of appropriate harmonic coordinates. In carrying out this program we prove moment and conductance estimates for MGW trees, which may be of independent interest. In addition, we extend our construction of the reversing measure to a biased random walk with random environment (RWRE) on MGW trees, again at a critical value of the bias. We compare this result against a transience--recurrence criterion for the RWRE generalizing a result of Faraud (2011) for Galton--Watson trees.
\end{abstract}

\maketitle

\section{Introduction}
\label{s:intro}

Let $\tree$ denote an infinite tree with root $\rt$. The $\lm$-biased random walk on $\tree$, hereafter denoted $\rw_\lm(\tree)$, is the Markov chain $(X_t)_{t\ge0}$ with $X_0=\rt$ such that given $X_t=v$ with offspring number $d_v$ and $v\ne\rt$, $X_{t+1}$ equals the parent of $v$ with probability $\lm/(\lm+d_v)$, and is uniformly distributed among the offspring of $v$ otherwise (and if $X_t=\rt$, then $X_{t+1}$ is uniformly distributed among the offspring of $\rt$).

For supercritical Galton--Watson trees without leaves, if $\rho$ denotes the mean offspring number, then $\rw_\lm$ is a.s.\ recurrent if and only if $\lm\ge\rho$ (\cite[Thm.~4.3 and Propn.~6.4]{\Lyons}), and ergodic if and only if $\lm>\rho$ (\cite[Propn.~9-131]{\KSK} and \cite[p.~944 and p.~954]{\Lyons}). With $\abs{v}$ denoting the (graph) distance from vertex $v$ to the root $\rt$, $\abs{X_t}/t$ converges a.s.\ to a speed $V$, with $V=V(\lm)$ deterministic, positive for $\lm<\rho$ and zero otherwise (see \cite{\LPPErgodic,\LPPBiased} for $\lm<\rho$ and \cite{\PZ} for $\lm=\rho$; the case $\lm>\rho$ follows trivially from positive recurrence).

Further, subject to no leaves and finite exponential moments for the offspring distribution, a quenched CLT for $\rw_\lm$ ($\lm\le\rho$) on single-type Galton--Watson trees was shown by Peres--Zeitouni \cite{\PZ}, and extended to the setting of random walk with random environment (RWRE) by Faraud \cite{\Faraud}. In contrast, if leaves occur, there emerges a zero-speed transient regime $\lm<\lm_c$ (for $\lm_c<\rho$) \cite{\LPPBiased} where the leaves ``trap'' the random walk and create slow-down. It follows from the results of Ben Arous et al.\ \cite{\GWLeaves} that in this setting, for sufficiently small $\lm$ there cannot be a (functional) CLT with diffusive scaling. Analogous results on the critical ($\rho=1$) Galton--Watson tree conditioned to survive were shown by Croydon et al.\ \cite{\CroydonFK}. In this paper we consider the critical case $\lm=\rho$, where \cite[Thm.~1]{\PZ} proves that on a.e.\ Galton--Watson tree, the processes $(\abs{X_{\flr{nt}}}/\sqrt{n})_{t\ge0}$ converge in law to the absolute value of a (deterministically) scaled Brownian motion. Their proof is based on the construction of harmonic coordinates and an explicit description of a reversing probability measure $\igwr$ for $\rw_\rho$ ``from the point of view of the particle.'' Having such an explicit description is a very delicate property: even for Galton--Watson trees, no such description is known for $\lm<\rho$ except at $\lm=1$ which is done by \cite[Thm.~3.1]{\LPPErgodic}. One thus might be led to believe that \cite[Thm.~1]{\PZ} is a particular property resulting from the independence inherent in the Galton--Watson law.

Here we show to the contrary that such a quenched CLT extends to the much larger family of supercritical multi-type Galton--Watson trees with finite type space. We allow for leaves (but condition on non-extinction), demonstrating that at $\lm=\rho$ the ``trapping'' phenomenon of \cite{\GWLeaves} does not arise. We also replace the assumption of exponential moments for the offspring distribution by an assumption of finite moments of order $\mmt>4$, so that our result restricted to the single-type case strengthens \cite[Thm.~1]{\PZ}. However, the main interest of our result lies in moving from an i.i.d.\ to a Markovian structure for the random tree.

As in \cite{\PZ}, the key ingredient in our proof is the construction of an explicit reversing (probability) measure $\imgwr$ for $\rw_\lm$ from the point of view of the particle, generalizing $\igwr$ to the multi-type setting, for $\lm$ at the critical value on the boundary between transience and recurrence. See \S\ref{s:imgwr} for the details of the construction which may be of independent interest.

The model we consider is as follows: let $\trees$ be the space of rooted trees with type, where each vertex $v$ is given a type $\type_v$ from a finite alphabet $\cQ$. We let $\cB_\trees$ be the $\si$-algebra on $\trees$ generated by the cylinder sets (determined by the restrictions of trees to finite neighborhoods of the root). We write $\tree$ for a generic element of $\trees$ and $\rt$ for its root. A \emph{multi-type Galton--Watson} tree is a random element $\tree\in\trees$, generated from a starting type $\type_\rt\in\cQ$  and a collection of probability measures $\q^a$ ($a\in\cQ$) on $$\cQ^\star \equiv \bigcup_{\ell\ge0} \cQ^\ell,$$ as follows: begin with a root vertex $\rt$ of type $\type_\rt$. Supposing inductively that the first $n$ levels of $\tree$ have been constructed, each vertex $v$ at the $n$-th level generates random offspring according to law $\q^{\type_v}$. For our purposes the ordering of the children does not matter, so each $\q^a$ may equivalently be regarded as a probability measure on configurations $\vec x = (x_b)_{b\in\cQ} \in(\Z_{\ge0})^\cQ$, where $x_b$ is the number of children of type $b$. Continuing to construct successive generations in this Markovian fashion, we denote the resulting law on $(\trees,\cB_\trees)$ by $\mgwtype{\type_\rt}$. We denote by $\mgw$ any mixture of the measures $(\mgwtype{a})_{a\in\cQ}$ (with \eqref{e:mgw.init} the canonical mixture) and let $\nx\equiv\set{\abs{\tree}<\infty}$ denote the event of extinction.

For $a,b\in\cQ$ let
$$A(a,b) = 
\sum_{\vec x} \q^a(\vec x) x_b,$$
the expected number of offspring of type $b$ at a vertex of type $a$. (Unless otherwise specified, the implicit assumption hereafter is that $\E_{\q^a}[\abs{\vec x}]<\infty$ for all $a\in\cQ$ where $\abs{\vec x}\equiv \sum_b x_b$.) Throughout the paper we will refer to the following assumptions:
\bnm
\item[\hypirred] The matrix $A\equiv (A(a,b))_{a,b\in\cQ}$ is irreducible with Perron--Frobenius eigenvalue $\rho$.
\item[\hypks]
$A$ is positive regular (every entry of $A^{n_0}$ is positive for some $n_0\in\N$), $\rho>1$, and $\E_{\q^a}[\abs{\vec x} \log \abs{\vec x}]<\infty$ for all $a\in\cQ$.
\item[\hypmmt{\mmt}] $\E_{\q^a}[\abs{\vec x}^\mmt]<\infty$ for all $a\in\cQ$.
\enm
Note that \hypirred\ and $\rho>1$ together imply $\mgwtype{a}(\nx)<1$ for all $a\in\cQ$.

\subsection{Central limit theorems}

We take all real-valued processes to be in the space $D[0,\infty)$ equipped with the topology of uniform convergence on compact intervals. Our main theorem is the following:

\bthm
\label{t:mclt}
Under \hypirred, \hypks, and \hypmmt{\mmt} with $\mmt>4$, for $\mgw$-a.e.\ $\tree\notin\nx$, if $X\sim\rw_\rho(\tree)$ then the processes $(\abs{X_{\flr{nt}}}/(\si\sqrt{n}))_{t\ge0}$ converge in law in $D[0,\infty)$ to the absolute value of a standard Brownian motion for $\si$ a deterministic positive constant (see \eqref{e:eta.sigma}).
\ethm

\brmk\label{r:polyg}
By \cite[Propn.~3.10.4]{\EK}, an equivalent statement is that the polygonal interpolation of $k/n\mapsto \abs{X_k}/(\si\sqrt{n})$ converges to standard Brownian motion in the space $C[0,\infty)$ (again with the topology of local uniform convergence).
\ermk

Let $\rw_\lm^\cts(\tree)$ denote the continuous-time version of $\rw_\lm(\tree)$, which when at $v\in\tree$ moves to the parent of $v$ (if $v\ne\rt$) at rate $\lm$ and to each offspring of $v$ at rate $1$.

\bcor\label{c:mgw.clt.cts}
Under the assumptions of Thm.~\ref{t:mclt}, for $\mgw$-a.e.\ $\tree\notin\nx$, if $X^\cts\sim\rw_\rho^\cts(\tree)$ then the processes
$(\abs{X^\cts_{nt}}/(\si\sqrt{2\rho n}))_{t\ge0}$ converge in law in $D[0,\infty)$ to the absolute value of a standard Brownian motion.
\ecor

\begin{figure}[h]
\caption{$\imgw_0$ tree}
\label{f:imgw}
\centering
\includegraphics[trim = 5cm 12cm 4cm 4cm,clip]{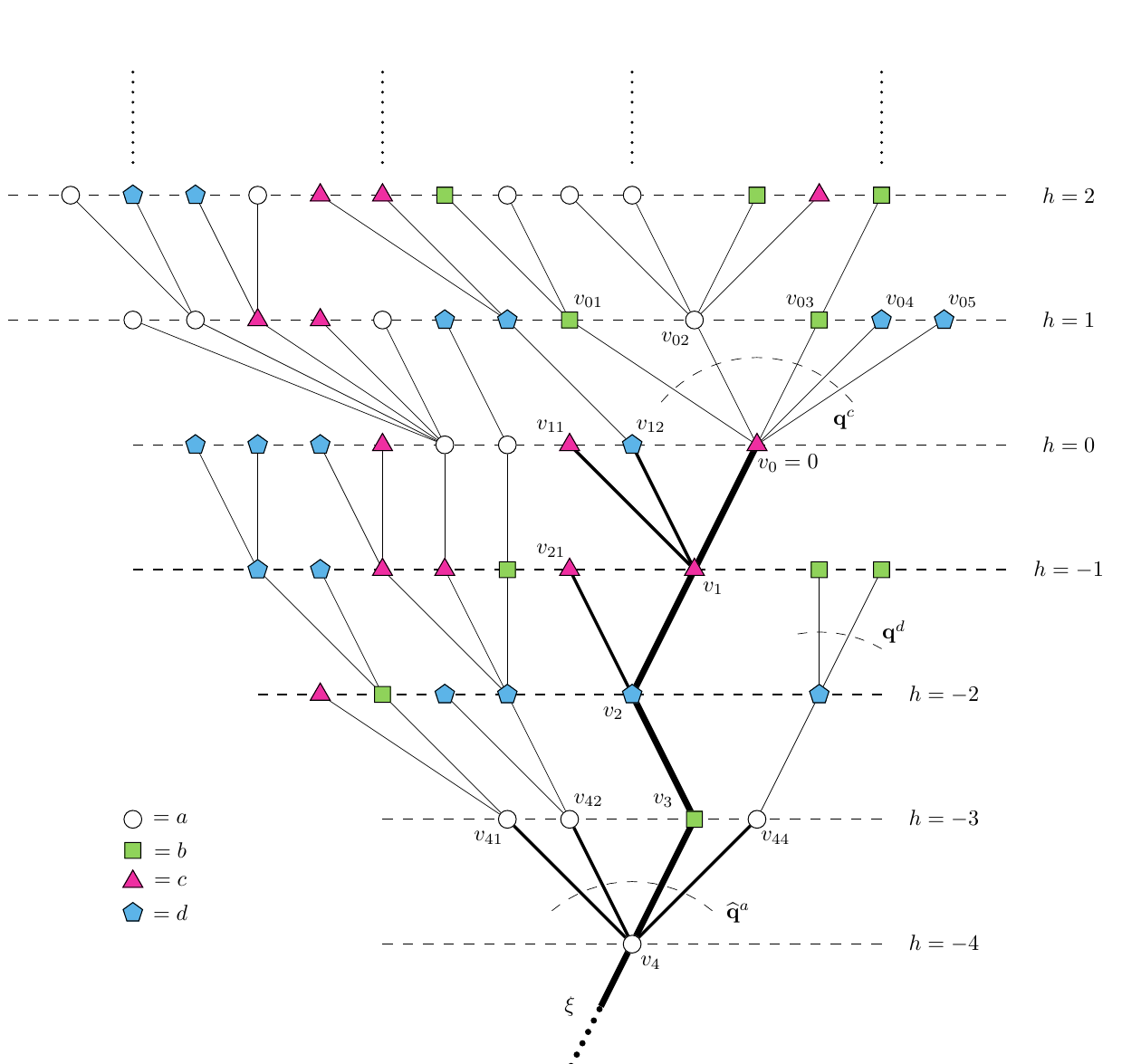}
\end{figure}

By moving the root of the tree to the current position of the random walk, $\rw_\lm$ on the tree induces a random walk on the space $\trees$, the ``walk from the point of view of the particle.'' As in \cite[\S3]{\PZ}, to make the latter process Markovian we amend the state space so as to keep track of the ancestry of the vertices. Specifically, we consider the space $\treesray$ of pairs $(\tree,\ray)$, where $\tree$ is an infinite tree and $\ray=(\rt=v_0,v_1,v_2,\ldots)$ is a ray emanating from the root $\rt$; this ray indicates the ancestry of each vertex in the tree. Let $\cB_\treesray$ denote the $\si$-algebra generated by the cylinder sets. We define a height function $h$ on $\tree$ as follows: set $h(v_n)=-n$, and for $v\notin\ray$ set
\beq \label{e:horocycle}
h(v) = h(R_v) + d(v,\ray)
\eeq
where $d$ denotes graph distance and $R_v$ is the nearest vertex to $v$ on $\ray$ (see Fig.~\ref{f:imgw}). We denote by $\rw_\lm(\tree,\ray)$ the $\lm$-biased random walk $(Y_t)_{t\ge0}$ on $(\tree,\ray)$, where the bias goes in the direction of decreasing height. With $\tree^v$ the tree $\tree$ rooted at $v$ instead of $\rt$, and $\ray^v$ the unique ray emanating from $v$ such that $\ray \cap \ray^v$ is an infinite ray, let
\beq \label{e:imgw.move.root}
(\tree,\ray)^{Y_t}\equiv(\tree^{Y_t},\ray^{Y_t}),\quad t\ge0.
\eeq
This is a Markov process with state space $\treesray$, and we hereafter refer to it as $\trw_\lm$. Let $\rw_\lm^\cts$ denote the continuous-time version of $\rw_\lm(\tree,\ray)$ (moving in the direction of increasing height at rate $1$ and in the direction of decreasing height at rate $\lm$), and let $\trw_\lm^\cts$ denote the induced continuous-time process on the space $\treesray$.

As in the single-type Galton--Watson case considered in \cite{\PZ}, the key to our proof lies in finding an explicit reversing measure $\imgw$ for $\trw_\rho^\cts$, which is then easily translated to a reversing measure $\imgwr$ for $\trw_\rho$. For a tree $\tree$ (with or without marked ray) and for any vertex $v\in\tree$, we denote by $\subtree{\tree}{v}$ the subtree induced by $v$ and its descendants, where descent is in direction of increasing distance from the root for a rooted tree, and in the direction of increasing height for a tree with marked ray. If $\mu$ is a law on trees we use $\mu\otimes\rw_\lm$ to denote the joint law of the tree together with the realization of $\rw_\lm$ on that tree.

\bthm \label{t:imgw} Assume \hypirred.
\bnm[(a)]
\item\label{t:imgw.rev}
There exists a reversing probability measure $\imgw$ for $\trw^\cts_\rho$, and if we define 
$$\f{d\imgwr}{d\imgw} = \f{d_\rt+\rho}{2\rho},$$
then $\imgwr$ is a reversing probability measure for $\trw_\rho$.
\item\label{t:imgw.erg}
For $((\tree,\ray),(Y_t)_{t\ge0}) \sim \imgwr\otimes\rw_\rho$, the stationary sequence $((\tree,\ray)^{Y_t})_{t\ge0}$ is ergodic.
\enm
\ethm

The $\imgw$ trees always have an infinite ray $\ray$, though the trees coming off the ray may be finite. The measures $\imgw,\imgwr$ are the multi-type analogues of the measures $\igw,\igwr$ of \cite{\PZ}. Thm.~\ref{t:imgw} and the construction of harmonic coordinates allow us to prove the following quenched CLT for $\rw_\rho$ on $\imgwr$ trees, which will be used to deduce Thm.~\ref{t:mclt}.

\bthm \label{t:iclt}
Under \hypirred, \hypks, and \hypmmt{\mmt} with $\mmt>2$, for $\imgwr$-a.e.\ $(\tree,\ray)$, if $Y\sim\rw_\rho(\tree,\ray)$ then the processes $(h(Y_{\flr{nt}})/(\si\sqrt{n}))_{t\ge0}$ converge in law in $D[0,\infty)$ to a standard Brownian motion.
\ethm

\subsection{Transience--recurrence boundary in random environment}

In the setting of $\rw_\lm$ on $\mgw$ trees, $\lm=\rho$ represents the onset of recurrence. Indeed, $\mgw$-a.e.\ tree $\tree$ on the event of non-extinction has branching number $\br\tree=\rho$ \cite[Propn.~6.5]{\Lyons}, therefore $\rw_\lm(\tree)$ is transient for $\lm<\rho$ and recurrent for $\lm>\rho$ \cite[Thm.~4.3]{\Lyons}. In fact, recurrence for all $\lm\ge\rho$ follows from a simple conductance calculation (for the general theory see \cite[Ch.~2]{lyons.peres}), therefore $\rho$ is the boundary between transience and recurrence for $\rw_\lm$ on $\mgw$ trees. Further $\rho$ is the boundary between non-ergodicity and ergodicity, with $\rw_\rho$ null recurrent (see \cite[Propn.~9-131]{\KSK} and \cite[p.~944 and p.~954]{\Lyons}) and of zero speed (e.g.\ from the bound of Lem.~\ref{l:carne}).

We believe that the existence of a reversing measure and CLT is a feature of the onset of recurrence in a more general setting. Indeed, suppose each vertex $v\in\tree\setminus\set{\rt}$ has, in addition to its type $\type_v$ from the (finite) alphabet $\cQ$, a weight $\al_v\in(0,\infty)$. Fixing such a tree $\tree$ (the environment), the $\lm$-biased random walk with random environment $\rwre_\lm(\tree)$ is the Markov chain $(X_t)_{t\ge0}$ with $X_0=\rt$ which, when at vertex $v$ with offspring
$$(\vec y,\vec \al)\equiv((y_1,\al_1),\ldots,(y_\ell,\al_\ell))\in\cQre^\ell,$$
jumps to a random neighbor $w$ of $v$ with probability proportional to $\al_w$ if $w$ is a child of $v$, and to $\lm$ if $w$ is the parent of $v$. (Note that $\rw_\rho(\tree)$ corresponds to the case $\al_v=1$ for all $v$.) We let $\rwre_\lm^\cts(\tree)$ denote the continuous-time version of $\rwre_\lm(\tree)$.

If $\q^a$ ($a\in\cQ$) is a probability measure on
$$\cQre^\star
\equiv \bigcup_{\ell\ge0} \cQre^\ell,\quad
\cQre = \cQ\times(0,\infty),$$
then the collection $(\q^a)_{a\in\cQ}$ together with starting type $\type_\rt\in\cQ$ specifies a law $\mgwretype{a_0}$ on the space $\trees$ of typed weighted rooted trees. As before we let $\mgwre$ denote any mixture of the $\mgwretype{a_0}$. This model, studied in the single-type case in \cite{\Faraud}, allows for quite general distributions on the (immediate) neighborhood of each vertex, but conditioned on types the weights in different neighborhoods must be independent.

For $\gam\in\R$ and $a,b\in\cQ$, let
\beq \label{e:rwre.mat}
\bar A^{(\gam)}(a,b)
\equiv \int_{\cQre^\star}
	\sum_j \Ind{y_j=b} \al_j^\gam
	\d\q^a(\vec y,\vec\al)
\eeq
(not necessarily finite for all $\gam$). Let $\bar\rho(\gam)$ be the Perron--Frobenius eigenvalue of $\bar A^{(\gam)}$ where well-defined (i.e.\ where $\bar A^{(\gam)}$ has finite entries and is irreducible), and $\infty$ otherwise. We will prove the following characterization of the transience--recurrence boundary for $\rwre_\lm$, extending part of \cite[Thm.~1.1]{\Faraud}:

\bthm\label{t:tr}
Suppose $\bar A^{(0)}$ is positive regular, and $\bar\rho(\gam)<\infty$ for $\gam$ in an open neighborhood of $0$. For $\lm>0$ let
$$p_\lm \equiv \inf_{0\le \gam\le 1} \f{\bar\rho(\gam)}{\lm^\gam}.$$
\bnm[(a)]
\item\label{t:tr.r} If $p_\lm<1$, then $\rwre_\lm$ is positive recurrent $\mgwre$-a.s.
\item\label{t:tr.t} If $p_\lm>1$, then $\rwre_\lm$ is transient $\mgwre(\cdot\giv\nx^c)$-a.s.
\enm
\ethm

Thus the transience--recurrence boundary for $\rwre_\lm$ occurs at the unique value $\lm=\rho^\circ$ for which $p_{\rho^\circ}=1$. On the other hand, let $\treesray$ denote the space of typed weighted trees with ray, and let $\trwre_\lm$ and $\trwre^\cts_\lm$ denote the Markov chains in $\treesray$ induced by $\rwre_\lm$ and $\rwre^\cts_\lm$ respectively. We have the following generalization of Thm.~\ref{t:imgw}~(\ref{t:imgw.rev}):
\bthm \label{t:rwre.imgw}
Suppose $\mgwre$ is such that $\bar A\equiv \bar A^{(1)}$ is irreducible with Perron--Frobenius eigenvalue $\bar\rho\equiv\bar\rho(1)$. Then there exists a reversing probability measure $\imgwre$ on $\treesray$ for $\trwre^\cts_{\bar\rho}$. If we let $\al_{0j}$ denote the weight for the $j$-th child of the root $\rt$, and set
$$\f{d\imgwrre}{d\imgwre} = \f{\bar\rho + \sum_{j=1}^{d_\rt} \al_{0j}}{2\bar\rho},$$
then $\imgwrre$ is a reversing probability measure for $\trwre_{\bar\rho}$.
\ethm

We can see that $\rho^\circ$ matches $\bar\rho$ if and only if the function $\gam\mapsto \bar\rho(\gam)/(\rho^\circ)^\gam$ attains its infimum over $0\le\gam\le1$ at $\gam=1$. If this fails, Thm.~\ref{t:rwre.imgw} still gives a reversing measure at $\bar\rho$, but $\bar\rho>\rho^\circ$ and the walk is already positive recurrent above $\rho^\circ$. However, at least in the single-type case, we have $\rho^\circ=\bar\rho$ in all cases in which a CLT is possible: indeed, if
$$\ka\equiv\inf\setb{\gam\ge0 : \f{\bar\rho(\gam)}{(\rho^\circ)^\gam}=1},$$
by results of \cite{\HuShi} a CLT cannot hold unless $\ka\ge2$ (see \cite[p.~3]{\Faraud}). We expect $\ka\ge2$ also to be a necessary condition in the multi-type case, and thus Thm.~\ref{t:tr} and Thm.~\ref{t:rwre.imgw} support the claim that reversing measures occur at the boundary between transience and recurrence \emph{in cases in which a CLT is possible}. However, even in the single-type case the random environment creates technical difficulties, and the RWRE-CLT of \cite{\Faraud} requires some restriction on $\ka$. While we expect that the methods of this paper and \cite{\Faraud} can also be adapted to extend the RWRE-CLT to the multi-type setting under the same restrictions on $\ka$, new ideas are required to achieve a CLT for the entire regime $\ka\ge2$.

\subsection*{Outline of the paper}

\bitm
\item In \S\ref{s:imgwr} we construct the reversing measure $\imgwr$ for $\trw_\rho$ (in \S\ref{ss:imgwr}) and its generalization $\imgwrre$ for $\trwre_{\bar\rho}$ (in \S\ref{ss:imgwr.re}); these constructions are based on ideas from \cite{\KLPP}. In \S\ref{ss:imgwr.wk.lim} we give an alternative characterization of $\imgwr$ (extending a characterization of \cite{\PZ} to the multi-type setting) which we use to prove ergodicity of the stationary sequence $((\tree,\ray)^{Y_t})_{t\ge0}$.

\item In \S\ref{s:qiclt} we prove the quenched $\imgwr$-CLT Thm.~\ref{t:iclt}: in \S\ref{ss:qiclt.harm} we construct on $\imgwr$-a.e.\ $(\tree,\ray)$ a function $v\mapsto S_v$ ($v\in\tree$) which is harmonic with respect to the transition probabilities of $\rw_\rho(\tree,\ray)$. By stationary and ergodicity of $((\tree,\ray)^{Y_t})_{t\ge0}$ with respect to $\imgwr$ we are able to control the quadratic variation of the martingale $M_t\equiv S_{Y_t}$ to obtain an $\imgwr$-a.s.\ martingale CLT. In \S\ref{ss:qiclt} we adapt the methods of \cite{\PZ} and \cite{\Faraud} to show that $h(Y_t)$ is uniformly well approximated by $M_t/\eta$ ($\eta$ an explicit constant), proving Thm.~\ref{t:iclt}.

\item In \S\ref{s:coup} we prove the quenched $\mgw$-CLT Thm.~\ref{t:mclt}. In \S\ref{ss:coup} we review (a slight modification of) a construction of \cite{\PZ} which gives a ``shifted coupling'' of $(\tree,(X_t)_{t\ge0})\sim\mgw\otimes\rw_\rho$ with $((\tree^\coup,\ray),(Y_t)_{t\ge0})\sim\imgw_0\otimes\rw_\rho$ such that fresh excursions of $X$ are matched with fresh excursions of $Y$ away from $\ray$. From this we obtain an \emph{annealed} $\mgw$-CLT (in \S\ref{ss:coup.amclt}) for $X$ by controlling the amount of time spent outside the coupled excursions as well as the drift of $Y$ along $\ray$. Because of the dependence between $\tree$ and $Y$ we do not see how to this coupling directly to prove a quenched ($\mgw$-a.s.) CLT. Instead, in \S\ref{ss:coup.qmclt} we adapt the method of \cite{\BolthausenSznitman} to deduce Thm.~\ref{t:mclt} from the annealed CLT by controlling the correlation between two realizations of $\rw_\rho$ on a single $\mgw$ tree $\tree$ (as was done in \cite[\S7]{\PZ} in the case $\lm<\rho$).

\item In \S\ref{s:rwre.trans.recur} we prove Thm.~\ref{t:tr} describing the transience--recurrence boundary for $\rwre_\lm$. The main result needed is a large deviations estimate (Lem.~\ref{l:rwre.large.devs}) on the conductances at the $n$-th level of the tree.

\item In \S\ref{s:appx} are collected some basic properties of $\mgw$ which are needed in the course of our proof and which may be of independent interest. In \S\ref{ss:appx.mmt} we show that moments for the offspring distributions translate directly to moments for the normalized population size defined in \S\ref{ss:imgwr.wk.lim}. In \S\ref{ss:appx.cond} we prove the existence of harmonic moments for the normalized population size, and use this result to prove conductance estimates used in the proof of Thm.~\ref{t:mclt}.
\eitm

\subsection*{Open problems}
We conclude this section by mentioning some open problems in this area. These problems are open even for single-type Galton--Watson trees.
\bnm
\item Does a CLT with diffusive scaling hold for $\rw_\rho$ in the entire regime $\mmt\ge2$?
\item Does a CLT with diffusive scaling hold for $\rwre_{\bar\rho}$ in the entire regime $\ka\ge2$?
\item What happens for simple random walk on the critical Galton--Watson tree (conditioned to survive)?
\item Does a CLT with any scaling (or other limit law) hold for $\rw_\rho$ when $\mmt<2$?
\enm
A common feature of these problems is that while the reversing measure for the process from the perspective of the particle is given by Thm.~\ref{t:imgw}, the method of martingale approximation used in \cite{\PZ,\Faraud} and in this paper seem not to be directly applicable.

\subsection*{Acknowledgements}
We are very grateful to Ofer Zeitouni for many helpful communications leading to the proof of ergodicity in Thm.~\ref{t:imgw} and the method of going from the annealed to quenched CLT in the proof of Thm.~\ref{t:mclt}. A.D.\ thanks Alexander Fribergh for discussions on the results of \cite{\GWLeaves} which motivated us to extend our results to trees with leaves. N.S.\ thanks Yuval Peres for several helpful conversations about the papers \cite{\PZ,\KLPP} which led to the construction of the reversing measure. We thank the anonymous referee for many valuable comments on drafts of this paper.

\section{Reversing probability measures for \texorpdfstring{$\trw_\rho$}{TRW(rho)} and \texorpdfstring{$\trwre_{\bar\rho}$}{TRWRE(rho)}}
\label{s:imgwr}

Assuming only \hypirred, in this section we construct the reversing measure $\imgwr$ for $\trw_\rho$ (\S\ref{ss:imgwr}) as well as its generalization $\imgwrre$ for $\trwre_{\bar\rho}$ (in \S\ref{ss:imgwr.re}). In \S\ref{ss:imgwr.wk.lim} we give an alternative characterization of $\imgwr$ which we use to prove ergodicity of the stationary sequence $((\tree,\ray)^{Y_t})_{t\ge0}$. Except in \S\ref{ss:imgwr.re} we work throughout with unweighted trees.

Consider a multi-type Galton--Watson measure $\mgw$ with offspring distributions $(\q^a)_{a\in\cQ}$ and mean matrix $A$. Hereafter we let $\rv$ and $\lv$ denote the right and left eigenvectors respectively associated to the Perron--Frobenius eigenvalue of $A$, normalized so that $\sum_a \lvec{a}=\sum_a \rvec{a}=1$. Since our results are stated for $\mgw$-a.e.\ tree, with no loss of generality we set hereafter
\beq\label{e:mgw.init}
\g(a)\equiv \mgw(\type_\rt=a)=\lvec{a}.
\eeq
Unless otherwise specified, $X$ and $Y$ denote $\rw_\rho$ on trees without and with marked ray respectively.

\subsection{Construction of \texorpdfstring{$\imgwr$}{IMGWR}}
\label{ss:imgwr}

We begin by constructing two auxiliary measures on the space $\treesray$ of trees with ray $(\tree,\ray)$. Let the infinite ray $\ray$ (without types) be given. For some $n>0$, we let vertex $v_n$ be given a type $\type_n$ according to a distribution $\bpi$, to be determined shortly. It is then given offspring $\vec x^{v_n}$ according to the \emph{inflated} offspring distribution
$\infq^{\type_n}$, where
$$\infq^a(\vec x)
\equiv \q^a(\vec x) 
\f{\ip{\vec x}{\rv}}
	{\rho \rvec{a}} \quad\forall a\in\cQ;$$
note that $\infq^a(\abs{\vec x}\ge1)=1$. One offspring $w$ of $v_n$ is then identified with the next vertex $v_{n-1}$ along $\ray$, where each $w$ is chosen with probability $\rvec{\type_w} / \ip{\vec x^{v_n}}{\rv}$. We proceed in this manner along the ray ending with the identification of $v_0=\rt$. The sequence of types $\type_n,\type_{n-1},\ldots$ seen along the ray
is then (by \hypirred) an irreducible Markov chain with transition probabilities
\beq \label{e:markov.trans}
K(a,b) = \sum_{\vec x} \infq^a(\vec x)
	\f{e_b x_b}
	{\ip{\vec x}{\rv}}
= \sum_{\vec x} \q^a(\vec x) \f{e_b x_b} {\rho \rvec{a}}
= \f{\rvec{b}}{\rho \rvec{a}} A(a,b).
\eeq
This chain has stationary distribution $\bpi(a)\equiv\rvec{a} \lvec{a}/\ip{\rv}{\lv}$, so starting with $\type_n \sim\bpi$ yields a consistent family of distributions for $(v_n,\ldots,v_1)$ and their (immediate) offspring, with types. By Kolmogorov's existence theorem, this uniquely specifies the distribution of the \emph{backbone} of the tree, that is, of the ray $\ray$ together with all (immediate) offspring of the vertices $v_i$, $i>0$. To each of these offspring (off the ray) and to $\rt$, we attach an independently chosen $\mgw$ tree conditioned on the given type, and denote by $\imgw_0$ the resulting measure on $\treesray$.

The \emph{inflated multi-type Galton--Watson} measure
$\imgw$ is obtained from $\imgw_0$ by an additional biasing
according to the root type $\type_\rt$. Specifically, we set
$$\f{d\imgw}{d\imgw_0}
= \f{1/\rvec{\type_\rt}}{\E_\bpi[1/\rvec{\type}]}
= \f{\E_\g[\rvec{\type}]}{\rvec{\type_\rt}},$$
where $\chi$ denotes a random variable on $\cQ$ with the specified distribution. We note that under $\imgw$, $\type_\rt\sim\g$ and so $\subtree{\tree}{\rt}$ has marginal law $\mgw$, which implies
\begin{equation*}
\E_\imgw[d_\rt] = \E_\mgw[d_\rt]
= \sum_a \lvec{a} \sum_b A(a,b) = \rho.
\end{equation*}
With this in mind, we define the probability measure $\imgwr$ such that
\beq\label{e:imgwr}
\f{d\imgwr}{d\imgw} = \f{d_\rt+\rho}{\E_\imgw[d_\rt+\rho]}
= \f{d_\rt+\rho}{2\rho},
\eeq
and proceed to show that it is a reversing measure for $\trw_\rho$. From now on we adopt the notation that if $\mu$ is a law on trees $\tree$ (with or without marked ray) and $a\in\cQ$, $\mu^a$ refers to the law conditioned on $\type_\rt=a$.

\bpf[Proof of Thm.~\ref{t:imgw}~(\ref{t:imgw.rev})]
For the purposes of this proof we let $\trees$ and $\treesray$ be spaces of \emph{labelled} (or planar) trees (without and with marked ray, respectively), with corresponding Borel $\si$-algebras $\cB_\trees$ and $\cB_\treesray$. We extend $\mgw$, $\imgw_0$, etc.\ to be measures on these spaces by choosing an independent uniformly random ordering for the offspring of each vertex. For $(\tree,\ray)\in\treesray$ we use the shorthand $i$ for $v_i\in\ray$, and write $(i1,\ldots,id_i)$ for its ordered offspring (with $\vec x^i$ denoting the counts of offspring of $i\equiv v_i$ of each type).

Recalling the notation of \eqref{e:imgw.move.root}, let $\cS$ denote the map $(\tree,\ray)\mapsto(\tree,\ray)^1$. We will show that for $A,B\in\cB_\treesray$, 
\beq \label{e:revers}
\int_A \mathbf{p}((\tree,\ray),B) \d\imgwr(\tree,\ray)
= \int_B \mathbf{p}((\tree',\ray'),A) \d\imgwr(\tree',\ray')
\eeq
where $\mathbf{p}((\tree,\ray),B)$ denotes the transition kernel of the process $\trw_\rho$. This identity implies reversibility of $\trw_\rho$ on the space of labelled trees. Since this process projects to $\trw_\rho$ on the space of unlabelled trees, the reversibility of the latter follows.

For $(\tree,\ray)\sim\imgw_0$, let $\barimgwztype{a}$ denote the law of the subtree $\tree\setminus \subtree{\tree}{i-1}$ rooted at $i$ with marked ray $\ray^i$, conditioned on the event $\set{\type_{i-1}=a}$, for any $i\ge1$ (note that this law does not depend on $i$). Then
$$d\imgw_0(\tree,\ray)
= \bpi(\type_1)
	\f{\infq^{\type_1}(\vec x^1)}{d_1!}
	\f{\rvec{\type_0}}{\ip{\vec x^1}{\rv}}
	d\barimgwztype{\type_1}(\tree\setminus \subtree{\tree}{1},\ray^2)
	\prod_{j=1}^{d_1} d\mgwtype{\type_{1j}}(\subtree{\tree}{1j}).$$
Let $\Pinj$ denote the collection of $\cB_\treesray$-measurable sets on which $\cS$ is injective, and suppose $B\in\Pinj$. If $\mu$ is a measure on $\treesray$, $\cS^*_B\mu(\cdot)\equiv (\mu\circ\cS)(B\cap\cdot)$ is a well-defined measure on $\treesray$. Then
$$d \cS^*_B\imgw_0(\tree,\ray)
= \Ind{(\tree,\ray)\in B}
\bpi(\type_1) \f{\q^{\type_1}(\vec x^1)}{d_1!}
d\barimgwztype{\type_1}(\tree\setminus \subtree{\tree}{1},\ray^2)
\prod_{j=1}^{d_1} d\mgwtype{\type_{1j}}(\subtree{\tree}{1j}).$$
so
$$\f{d\cS^*_B\imgw_0}{d\imgw_0}
= \I_B \f{\rho \rvec{\type_1}}{\rvec{\type_0}}.$$
We then verify that
\begin{align}
\nonumber
\f{d\cS^*_B\imgw}{d\imgw}
&= \I_B
\f{ \big( \f{d\cS^*_B\imgw}{d\cS^*_B\imgw_0} \big)}
	{\big( \f{d\imgw}{d\imgw_0} \big)}
	\f{d \cS^*_B\imgw_0}{d\imgw_0}\\
\label{e:imgw.rev} &= \I_B
	\f{\big( \f{d\imgw}{d\imgw_0}\circ\cS \big)}
	{\big( \f{d\imgw}{d\imgw_0} \big)}
	\f{d \cS^*_B\imgw_0}{d\imgw_0}
= \I_B \f{1/\rvec{\type_1}}{1/\rvec{\type_\rt}} \f{\rho\rvec{\type_1}}{\rvec{\type_\rt}}
= \I_B\rho,
\end{align}
and similarly
$$\f{d\cS^*_B\imgwr}{d\imgwr}=\I_B\f{\rho(d_1+\rho)}{d_\rt+\rho}.$$
The left-hand side of \eqref{e:revers} can be written as
$$\int_{A\cap \cS^{-1} B} \f{\rho}{d_\rt+\rho} \d\imgwr(\tree,\ray)
+ \int_A \f{1}{d_\rt+\rho} \sum_{i=1}^{d_\rt}
	\Ind{(\tree^{0i},\ray^{0i})\in B}\d\imgwr(\tree,\ray).$$
Using the injectivity of $\cS$ on $B$, the second integral can be written as
\begin{align*}
\int_{A\cap \cS B} \f{1}{d_\rt+\rho} \d\imgwr(\tree,\ray)
&= \int_{\cS^{-1} A \cap B} \f{1}{d_1+\rho} \d\cS^*_B\imgwr(\tree,\ray) \\
&= \int_{\cS^{-1} A \cap B} \f{\rho}{d_\rt+\rho} \d\imgwr(\tree,\ray).
\end{align*}
Combining these yields an expression for the left-hand side of \eqref{e:revers} which is symmetric in $A$ and $B$, from which it is clear that the two sides must agree.

Since every cylinder event $F$ can be decomposed into the disjoint union of the event $F_j = \setb{(\tree,\ray)\in F : \rt = 1j}$ (i.e., $\rt$ is the $j$-th child of $1$), with $F_j$ clearly in $\Pinj$, we have that $\Pinj$ generates $\cB_\treesray$. To conclude, for fixed $A$ let $\cB_\treesray'$ denote the collection of sets $B\in\cB_\treesray$ for which \eqref{e:revers} holds. From the above $\cB_\treesray'$ contains the $\pi$-system $\Pinj$. Further $\cB_\treesray'$ is closed under monotone limits and countable disjoint unions, and in particular it contains $\treesray$ since $\treesray$ can be decomposed as a countable disjoint union of sets in $\Pinj$ by a similar argument as above. Thus by the $\pi$-$\lm$ theorem \eqref{e:revers} holds for all $B\in\si(\Pinj)$, and extends to all $B\in\cB_\treesray$ again using the claim above.

The proof that $\imgw$ is a reversing measure for the Markov pure jump process $\trw^\cts_\rho$ is similar: instead of \eqref{e:revers} we show that
\beq\label{e:revers.cts}
\int_A \lm(\tree,\ray) \mathbf{p}((\tree,\ray),B) \d\imgw(\tree,\ray)
= \int_B \lm(\tree,\ray) \mathbf{p}((\tree',\ray'),A) \d\imgw(\tree',\ray')
\eeq
where $\lm(\tree,\ray)\equiv\lm+d_\rt$ is the instantaneous jump rate of the process at state $(\tree,\ray)$. As before, it suffices to show this for $B\in\Pinj$. In this case the left-hand side of \eqref{e:revers.cts} equals
$$\int_{A\cap \cS^{-1} B} \rho \d\imgw(\tree,\ray)
+ \int_{\cS^{-1} A \cap B} \d\cS^*_B\imgw(\tree,\ray),$$
which by \eqref{e:imgw.rev} coincides with the right-hand side of \eqref{e:revers.cts}.
\epf

\subsection{Extension of \texorpdfstring{$\imgwr$}{IMGWR} to random environment}
\label{ss:imgwr.re}

We now extend the methods of the previous section to prove Thm.~\ref{t:rwre.imgw}. Let $\rvre,\lvre$ denote the right and left Perron--Frobenius eigenvectors of $\bar A \equiv \bar A^{(1)}$, normalized to have sum $1$; as before we set $\g(a)\equiv\mgwre(\type_\rt=a)$ to be $\lvecre{a}$.

We proceed much as in the deterministic environment setting, although the notation becomes more complicated. For $\vec y\in\cQ^\ell$ write $\rvre(\vec y)\equiv(\rvecre{y_j})_{j=1}^\ell$. For $a\in\cQ$
$$\E_{\q^a} [\ip{\rvre(\vec y)}{\vec \al}]
= \sum_b \bar A(a,b) \rvecre{b}
= \bar\rho \rvecre{a},$$
so we define the inflated offspring measure $\infq^a$ by
$$\f{d\infq^a}{d \q^a}
= \f{\ip{\rvre(\vec y)}{\vec\al}}{\bar\rho \rvecre{a}}.$$
We then construct the measure $\imgwre_0$ on $\treesray$ generalizing the measure $\imgw_0$ of the previous section: let the infinite ray $\ray$ (without types or weights) be given, and for some $n>0$ let $v_n$ have type $\type_n$. It is given offspring $(\vec y^{v_n},\vec\al^{v_n}) \sim \infq^{\type_n}$. One offspring $w$ of $v_n$ is identified with the next vertex $v_{n-1}$ along $\ray$, where each $w$ is chosen with probability
$$\f{\rvecre{\type_w} \al_w}{\ip{\rvre(\vec y)}{\vec\al^{v_n}}}.$$
Continuing the procedure along the ray up to $v_0=\rt$, the sequence of types $\type_n,\type_{n-1},\ldots$ seen along $\ray$ is an irreducible Markov chain with transition probabilities
$$\bar K(a,b)
= \E_{\infq^a}\bigg[ \f{\rvecre{b} \sum_j \al_j \Ind{y_j=b}}{\ip{\rvre(\vec y)}{\vec\al}} \bigg]
= \f{\rvecre{b}}{\bar\rho \rvecre{a}} \bar A(a,b)$$
and stationary distribution $\bar\bpi(a)=\rvecre{a}\lvecre{a}/\ip{\rvre}{\lvre}$. Thus, starting with $\type_n\sim\bar\bpi$ and applying Kolmogorov's existence theorem, we obtain a measure $\imgwre_0$ on $\treesray$ which is a generalization of $\imgw_0$.

\bpf[Proof of Thm.~\ref{t:rwre.imgw}]
The proof is by a straightforward modification of the proof of Thm.~\ref{t:imgw}~(\ref{t:imgw.rev}). Let $\cS:(\tree,\ray)\mapsto(\tree,\ray)^1$; we emphasize that $\cS$ is a mapping on typed \emph{weighted} labelled trees. For $(\tree,\ray)\sim\imgwre_0$, let $\barimgwztypere{a}$ denote the law of the subtree $\tree\setminus\subtree{\tree}{i-1}$ rooted at $i$ with marked ray $\ray^i$, conditioned on the event $\set{\type_{i-1}}$, for any $i\ge1$. Let $\Pinj$ denote the collection of $\cB_\treesray$-measurable sets on which $\cS$ is injective. For $A\in\cB_\treesray$ and $B\in\Pinj$, we compute
\begin{align*}
d\imgwre_0(\tree,\ray)
&= \bar\bpi(\type_1)
 \infq^{\type_1}(\vec y^1,\vec\al^1)
 \f{\rvecre{\type_0} \al_0}{\ip{\rvre(\vec y^1)}{\vec\al^1}}
 d\barimgwztypere{\type_1}(\tree\setminus \subtree{\tree}{1},\ray^2)
 \prod_{j=1}^{d_1} d\mgwretype{\type_{1j}}(\subtree{\tree}{1j}), \\
d \cS^*_B \imgwre_0(\tree,\ray)
&= \Ind{(\tree,\ray)\in B} \bar\bpi(\type_1)
	\q^{\type_1}(\vec y^1,\vec\al^1)
	d\barimgwztypere{\type_1}(\tree\setminus\subtree{\tree}{1},\ray^2)
 \prod_{j=1}^{d_1} d\mgwretype{\type_{1j}}(\subtree{\tree}{1j}),
\end{align*}
so $\rvecre{\type_0}\al_0 
\d\cS^*_B \imgwre_0 = \I_B 
\bar\rho \rvecre{\type_1} \d\imgwre_0$. Letting
\begin{align*}
\f{d\imgwre}{d\imgwre_0}
&\equiv \f{1/\rvecre{\type_\rt}}{\E_{\bar\bpi}[1/\rvecre{\type}]}
= \f{\E_{\ol\g}[\type]}{\rvecre{\type_0}},\\
\f{d\imgwrre}{d\imgwre}
&\equiv \f{\bar\rho + \sum_{j=1}^{d_\rt} \al_{0j}}
	{\E_{\imgwre}[\bar\rho + \sum_{j=1}^{d_\rt} \al_{0j}]}
= \f{\bar\rho + \sum_{j=1}^{d_\rt} \al_{0j}}
	{2\bar\rho},
\end{align*}
we obtain
\begin{align*}
\al_0 \d\cS^*_B\imgwre
&=\I_B\bar\rho\d\imgwre,\\
\f{\al_0}{\bar\rho + \sum_{j=1}^{d_1} \al_{1j}} \d\cS^*_B\imgwrre
&= \I_B \f{\bar\rho}{\bar\rho+\sum_{j=1}^{d_0}\al_{0j}} \d\imgwrre.
\end{align*}
The analogue of \eqref{e:revers} thus holds for all $B\in\Pinj$, and we extend to all $B\in\cB_\treesray$ by essentially the same argument used in the proof of Thm.~\ref{t:imgw}~(\ref{t:imgw.rev}).
\epf

\subsection{\texorpdfstring{$\imgw_0$}{IMGW(0)} as a weak limit and ergodicity}
\label{ss:imgwr.wk.lim}

In this section we provide an alternative characterization (Propn.~\ref{p:imgw.wk}) of the inflated Galton--Watson measure $\imgw_0$, which is then used in proving the ergodicity result Thm.~\ref{t:imgw}~(\ref{t:imgw.erg}). Propn.~\ref{p:imgw.wk} is also of independent interest as a multi-type extension of \cite[Lem.~1]{\PZ}.

To this end, we will define the notion of ``normalized population size'' for rooted trees $\tree$ with type. Let $\tree_n$ denote the subtree induced by $\setb{v\in\tree:\abs{v}\le n}$, and $D_n$ the set $\setb{v\in\tree:\abs{v}=n}$. Let $(\cF_n)_{n\ge0}$ denote the natural filtration of the tree, i.e., $\cF_n$ is the $\si$-algebra generated by $\tree_n$ (a finite tree with vertex types). Let $\vec Z_n = (Z_n(b))_{b\in\cQ}\in(\Z_{\ge0})^\cQ$ count the number of vertices of each type at level $n$, so $\vec Z_n$ is $\cF_n$-measurable. Then
$$\Znorm_n \equiv \f{\ip{\vec Z_n}{\rv}}{\rho^n}
= \f{1}{\rho^n} \sum_{v\in D_n} \rvec{\type_v}$$
is a non-negative $(\cF_n)$-martingale under $\mgwtype{a}$ for every $a$, with $\E_{\mgwtype{a}}[\Znorm_0]=\rvec{a}$ (see e.g.\ \cite[p.~49]{\Harris}). By the normalized population size of the tree we mean the a.s.\ limit of $\Znorm_n$, denoted  $W_\rt$. For $v\in \tree$ we use $W_v$ to denote the normalized population size of $\subtree{\tree}{v}$. Under \hypirred\ and \hypks, it follows from the multi-type Kesten--Stigum theorem (see \cite{\KSMGW}, or the conceptual proof of \cite{\KLPP}) that $W_\rt>0$ a.s.\ on the event of non-extinction, and $\E_{\mgwtype{a}}[W_\rt] = \rvec{a}$.

For $a\in\cQ$ let $\Qntype{a}_n$ be a probability measure on (infinite) rooted trees defined by
\beq\label{e:Qn.a}
\f{d\Qntype{a}_n}{d\mgwtype{a}} = \f{\Znorm_n}{\rvec{a}}.
\eeq
For $\tree\sim\Qn^a_n$ choose $v_n \in D_n$ at random with probabilities proportional to weights $\rvec{\type_{v_n}}$, and let $\Qn^a_{n\star}$ denote the law of the resulting pair $(\tree,v_n)$. Let $\Qn_{n\star} \equiv \sum_{a\in\cQ} \bpi_a \Qntype{a}_{n\star}$ and $\Qn_n \equiv \sum_{a\in\cQ} \bpi_a \Qntype{a}_n$, so that $d\Qn_n/d\mgw = \Znorm_n / \E_\g[\rvec{\type}]$. Finally let $\imgw_0(n)$ denote the law of $(\tree,\ray_0)^{v_n}$ (see \eqref{e:imgw.move.root} for this notation), where $(\tree,v_n)\sim\Qn_{n\star}$ and $\ray_0$ is any infinite ray emanating from $\rt$ not sharing an edge with the geodesic from $\rt$ to $v_n$.

\bppn\label{p:imgw.wk}
Under \hypirred, $\imgw_0(n)$ converges weakly to $\imgw_0$.
\eppn

The proposition can be seen from the following explicit construction of $\Qntype{a}_{n\star}$: begin with $v_0\equiv\rt$ of type $a$, and suppose inductively that we have constructed $(\tree_i,v_i)$ ($i<n$) where $\tree_i$ is the tree up to level $i$ and $v_i$ is the $i$-th vertex on the geodesic from $\rt$ to $v_n$. Then $v_i$ is given offspring $\vec x^{v_i}$ according to $\infq^{\type_{v_i}}$, and one of these offspring $w$ is randomly chosen (according to weights $\rvec{w}$) to be distinguished as $v_{i+1}$. Meanwhile all other vertices $v\in D_i \setminus \set{v_i}$ are given offspring $\vec x^v$ according to $\q^{\type_v}$. Once $(\tree_n,v_n)$ has been constructed, attach to each $v\in D_n$ an independent $\mgwtype{\type_v}$ tree. For $N\ge n$,
$$\f{\Qntype{a}_{n\star}(\tree_N,v_n)}{\mgwtype{a}(\tree_N)}
= \prod_{i=0}^{n-1} \f{\ip{\vec x^{v_i}}{\rv}}
	{\rho \rvec{\type_{v_i}}}
	\f{\rvec{\type_{v_{i+1}}}}{\ip{\vec x^{v_i}}{\rv}}
= \f{\rvec{\type_{v_n}}}{\rho^n \rvec{a}},$$
and summing over $v_n\in D_n$ gives \eqref{e:Qn.a}.

Letting $n\to\infty$ in $\Qntype{a}_{n\star},\Qn_{n\star}$ we obtain the measures $\Qntype{a}_{\infty\star},\Qn_{\infty\star}$ on rooted trees with infinite marked ray which coincide precisely with the measures $\infmgwstartype{a},\infmgwstar$ of \cite{\KLPP}. The corresponding marginals $\Qn^a_\infty\equiv\infmgwtype{a},\Qn_\infty\equiv\infmgw$ on trees without marked ray satisfy
$$\left.
	\f{d\Qn^a_\infty}
	{d\mgwtype{a}}\right|_{\cF_n}
= \f{\Znorm_n}{\rvec{a}},\quad
\left.\f{d\Qn_\infty}
	{d\mgw}\right|_{\cF_n}
= \f{\Znorm_n}{\E_\g[\rvec{\type}]}.$$
By the Kesten--Stigum theorem and Scheff\'e's lemma (see e.g.\ \cite[\S5.10]{\Williams}), $\Znorm_n \stackrel{L^1}{\longrightarrow} W_\rt$, hence
$$\f{d\Qn^a_\infty}{d\mgwtype{a}}
=\f{W_\rt}{\rvec{a}},\quad
\f{d\Qn_\infty}{d\mgw}
=\f{W_\rt}{\E_\mgw[W_\rt]}
=\f{W_\rt}{\E_\g[\rvec{\type}]}.$$
We remark that although $\Qn_{\infty\star}\equiv\infmgwstar$ and $\imgwr$ are both measures on trees with rays, they are not in general equivalent unless $K$ is reversible.

\bpf[Proof of Propn.~\ref{p:imgw.wk}]
Since $d\Qn_n/d\mgw = \Znorm_n / \E_\g[\rvec{\type}]$ and $\type_\rt\sim\g$ under $\mgw$, it follows that $\type_\rt\sim\bpi$ under $\Qn_{n\star}$. It is then clear from the constructions of $\imgw_0$ and $\Qn_{n\star}$ that if $(\tree,\ray)\sim\imgw_0$, then $(\subtree{\tree}{v_n},\rt)\sim\Qn_{n\star}$. In other words the portion of $(\tree,\ray)$ descended from $v_n$ has the same distribution under $\imgw_0(n)$ as under $\imgw_0$, proving the result.
\epf

Turning now to the proof of Thm.~\ref{t:imgw}~(\ref{t:imgw.erg}), it is useful to define a two-sided version of $\imgw_0$, as follows. Let $\treesline$ denote the space of trees with marked line: pairs $(\tree,\lne)$ where $\tree$ is an infinite tree and
$$\lne\equiv(\ldots,\lne_{-1},\lne_0=\rt,\lne_1,\ldots)$$
is a \emph{line} (doubly infinite simple path) passing through the root. The positive and negative parts $\lne_\pm\equiv(\lne_{\pm j})_{j\ge0}$ of $\lne$ are edge-disjoint rays emanating from $\rt$.

Now suppose in the construction of $\imgw_0$ we continue the backbone indefinitely rather than stopping at $\rt$, so that Kolmogorov's existence theorem gives a doubly infinite backbone based on a line $\lne$. Attaching $\mgw$ trees to the leaves of this backbone then gives a tree with marked line $(\tree,\lne)$, whose law $\imgwline$ is clearly stationary with respect to the shift $\ul\cS:(\tree,\lne)\mapsto(\tree,\lne)^{\lne_{-1}}$ which defined by moving the root to $\lne_{-1}$. (Alternatively, if $(\tree,\ray)$ has law $\imgw_0$ conditioned on non-extinction of $\subtree{\tree}{\rt}$ and $\lne$ is any line with $\lne_-=\ray$, then $\ul\cS^n(\tree,\lne)$ converges weakly to an $\imgwline$ tree.)

It follows from the discussion preceding Propn.~\ref{p:imgw.wk} that if we let
$$(\tree,\ray)\sim\imgw_0\quad\text{and}\quad
(\tree',\ray')\sim\Qn^{\type_\rt}_{\infty\star}$$
(independently conditioned on $\type_\rt$), and we delete from $\tree$ all the vertices descended from $\rt$ and identify $\rt$ with the root of $\tree'$, then we obtain a tree with marked line $\lne_-=\ray$, $\lne_+=\ray'$ whose law is precisely $\imgwline$. It follows that the marginal law $\imgw_1$ of $(\tree,\lne_-)$ under $\imgwline$ is given by
\beq\label{e:imgw1}\f{d\imgw_1}{d\imgw_0}
= \f{d\Qn^{\type_\rt}_\infty}{d\mgwtype{\type_\rt}}
= \f{W_\rt}{\rvec{\type_\rt}}.
\eeq

\bpf[Proof of Thm.~\ref{t:imgw}~(\ref{t:imgw.erg})]
We adapt the proof of \cite[Cor.~2.1.25]{\Zeitouni}. Abbreviating $\aT\equiv(\tree,\ray)$, we let $\nu$ denote the law of $\vec\aT\equiv(\aT_t)_{t\ge0}\equiv((\tree,\ray)^{Y_t})_{t\ge0}$ in the space $\treesrayinf$ of sequences of trees with ray, and $\shift_0$ the shift $(\aT_0,\aT_1,\ldots)\mapsto(\aT_1,\aT_2,\ldots)$ on $\treesrayinf$. The content of the result is that the measure-preserving system $(\treesrayinf,\cF^\infty,\nu,\shift_0)$ is ergodic.

\medskip\noindent
{\it Step 1: reduction to induced system.} \\
Recall that under the measure $\imgw_0$ the trees $\subtree{\tree}{i}\setminus\subtree{\tree}{i-1}$ are conditionally independent given the ray $\ray$ with types, and $\max_{a\in\cQ}\mgw^a(\nx)<1$; therefore it holds $\imgw_0$-a.s. that $\abs{\subtree{\tree}{i}}$ for infinitely many $i\in\ray$. Since the walk $Y_t$ on $(\tree,\ray)$ has a backward drift along $\ray$, this implies that if we let
$$A\equiv\set{\vec\aT\in\treesrayinf:
	\aT_0=(\tree,\ray) \text{ with }
	\abs{\subtree{\tree}{\rt}}=\infty}$$
and $n_A(\vec\aT)\equiv\inf\set{n\ge1 : \shift_0^n\vec\aT \in A}$ the first hitting time of $A$ after time zero, then $\nu(n_A<\infty)=1$. Thus $(\treesrayinf,\nu,\shift_0)$ forms a (Kakutani) tower over the induced measure-preserving system $(A,\nu_A\equiv\nu(\cdot\giv A),\shift_0^{n_A})$. We now show that the induced system is ergodic, which is equivalent to ergodicity of the original system (\cite{\Petersen}; see also \cite[\S2]{\LPPErgodic}).

\medskip\noindent
{\it Step 2: reduction to $\cS$-invariance.} \\
Let $n^i_A(\vec\aT)$ denote the $i$-th hitting time of $A$ after time zero and $\cH_i\equiv\si(\aT_0,\ldots,\aT_{n_A^i})$; note that $(\aT_{n_A^i})_{i\ge0}$ forms an $(\cH_i)$-Markov chain. Write $\shift\equiv\shift_0^{n_A}$ and let $\cI$ denote the $\si$-field of $\shift$-invariant subsets of $A$. Fix $B\in\cI$, and define
$$\phi:\treesray\to[0,1],\quad
\phi(\aT) \equiv \nu_A (
	\vec\aT\in B \giv \aT_0=\aT ).$$
The $\shift$-invariance of $B$ together with the Markov property implies
$$\nu_A[\vec\aT \in B\giv \cH_i]
= \nu_A[\shift^i\vec\aT \in B\giv \cH_i]
= \nu_A[\shift^i\vec\aT \in B\giv \aT_{n_A^i}]
= \phi(\aT_{n_A^i}),$$
i.e., $\phi(\aT_{n_A^i})$ is an $(\cH_i)$-martingale. By L\'evy's upward theorem, $\lim_{i\to\infty} \phi(\aT_{n_A^i}) = \I_B$, $\nu_A$-a.s., so that for any $0<a\le b<1$,
$$\f{1}{t}\sum_{i=0}^{t-1} \Ind{\phi(\aT_{n_A^i}) \in [a,b]}$$ converges $\nu_A$-a.s.\ to zero. On the other hand, by the Birkhoff ergodic theorem (see e.g.\ \cite[Thm.~6.2.1]{durrett}), it converges $\nu_A$-a.s.\ to $\nu_A (\phi(\tree)\in[a,b]\giv \cI)$. Taking expectations on both sides we find $\phi(\tree)\in\set{0,1}$ $\nu_A$-a.s., that is, $\phi=\I_{C_0}$ for some $C_0\in\cB_\treesray$. Further, since $\phi$ is a $\set{0,1}$-valued martingale, it holds $\nu$-a.s.\ that $\aT_0\in C_0$ if and only if $\aT_{n_A}\in C_0$. Since $\nu_A(\aT_{n_A}=\cS\aT_0)>0$ where $\cS$ is as defined above, $\I_{C_0}\le\I_{\cS^{-1} C_0}$ $\nu_A$-a.s. Applying the same argument with the martingale $1-\phi$ gives $\I_{C_0}=\I_{\cS^{-1} C_0}$ $\nu_A$-a.s., i.e., that $C_0$ is $\cS$-invariant.

\medskip\noindent
{\it Step 3: $\imgwr$-triviality of $\cS$-invariant sets.} \\ It follows from $B\subseteq A$ that $C_0$ is a subset of $A_0\equiv\set{(\tree,\ray)\in\treesray:\abs{\subtree{\tree}{\rt}}=\infty}$. Since $\imgwr\ll\imgw_1$ on $A_0$ by \eqref{e:imgw1}, the result follows by showing that $\cS$-invariant subsets of $A_0$ are $\imgw_1$-trivial. For any $C_0'\subseteq A_0$ let $\ul C_0'\equiv\set{(\tree,\lne):(\tree,\lne_-)\in C_0}$; the $\cS$-invariance of $C_0$ implies $\ul\cS$-invariance of $\ul C_0$. But the ergodicity of the Markov chain of types along the line $\lne$ readily implies that $\ul\cS$-invariant subsets of $\ul A_0$ are $\imgwline$-trivial, e.g.\ by the following modification of the argument of \cite[Thm.~2.15]{\Kallenberg}: take $\ul C_0^n$ measurable with respect to the portion $\lne_{[-n,0]}$ of the line between $\lne_{-n}$ and $\rt$, together with the descendant subtrees of $\lne_{-n},\ldots,\lne_{-1}$ away from $\lne$, such that the symmetric difference $\ul C_0^n \triangle \ul C_0$ has $\imgwline$-measure tending to zero in $n$. It follows from $\ul\cS$-invariance of $\ul C_0$ together with $\ul\cS$-stationarity of $\imgwline$ that $\imgwline[\ul C_0^n \triangle \ul C_0]=\imgwline[(\ul\cS^m\ul C_0^n) \triangle \ul C_0]$ for any $m$, so (by the triangle inequality)
$$\lim_{n\to\infty} \sup_m
	\absb{\imgwline(\ul C_0)
	-\imgwline[\ul C_0^n\cap(\ul\cS^m \ul C_0^n)]} = 0.$$
But for any $m>n$ we have
$$\imgwline[\ul C_0^n\cap(\ul\cS^m \ul C_0^n)]
= \imgwline[\lne_{[-n,0]},\lne_{[m-n,m]}]
\imgwline[\ul C_0^n\giv \lne_{[-n,0]}]
\imgwline[\ul\cS^m \ul C_0^n\giv \lne_{[m-n,m]}],$$
which tends as $m\to\infty$ to
$\imgwline[\ul C^n_0]^2$. Therefore
$$\imgwline[\ul C_0]
=\lim_{n\to\infty}\imgwline[\ul C^n_0]
=\lim_{n\to\infty}\imgwline[\ul C^n_0]^2
=\imgwline[\ul C_0]^2$$
which gives $\imgwline[\ul C_0]\in\set{0,1}$ as required.
\epf

\section{Harmonic coordinates and quenched \texorpdfstring{$\imgwr$}{IMGWR}-CLT}
\label{s:qiclt}

In this section we prove the quenched $\imgwr$-CLT Thm.~\ref{t:iclt}. Let
\beq \label{e:eta.sigma}
\eta\equiv\E_{\Qn_\infty}[W_\rt]
= \f{\E_\mgw[W_\rt^2]}{\E_\g[\rvec{\type}]},\quad
\si^2\equiv\f{\E_\g[\rvec{\type}]^2}{\E_\mgw[W_\rt^2]}.
\eeq
In \S\ref{ss:qiclt.harm} we construct harmonic coordinates for $\rw_\rho$ on $\imgwr$-a.e.\ $(\tree,\ray)$, and use the ergodicity result Thm.~\ref{t:imgw}~(\ref{t:imgw.erg}) proved above to show an $\imgwr$-a.s.\ CLT for the martingale $M_t\equiv S_{Y_t}$, with $M_{\flr{nt}}/(\eta\si\sqrt{n})$ converging to standard Brownian motion. In \S\ref{ss:qiclt} we control the error between $h(Y_t)$ and $M_t/\eta$ to prove Thm.~\ref{t:iclt}. The following result, whose proof is deferred to \S\ref{ss:appx.mmt}, implies finiteness of $\eta$ and $\si$ under \hypmmt{2}:

\bppn \label{p:moments}
If \hypirred, \hypks, and \hypmmt{\mmt} hold with $\mmt>1$, then $\E_\mgw[W_\rt^\mmt]<\infty$.
\eppn

\subsection{Harmonic coordinates for \texorpdfstring{$\rw_\rho$}{RW(rho)} and martingale CLT}
\label{ss:qiclt.harm}

From now on, if $\mu$ is a probability measure on trees (with or without marked ray), we use $\mu$ as shorthand also for $\mu\otimes\rw_\rho$. We write $\P_\tree$ for the law of the quenched random walk $\rw_\rho(\tree)$ and $\E_\tree$ for expectation with respect to $\P_\tree$, and let $(\cG^\tree_t)_{t\ge0}$ denote the corresponding filtration of the walk. Given $\tree$, for a vertex $v\in\tree$ we let $\pd v$ denote the neighbors of $v$, and $\desc v$ the offspring of $v$, i.e., $\desc v=\pd v\cap\subtree{\tree}{v}$. We write $v\le w$ if $w\in\subtree{\tree}{v}$, with $v<w$ if $w\ne v$.

For $v\in\tree$ recall that $W_v$ denotes the normalized population size of the subtree $\subtree{\tree}{v}$.\footnote{Note that if $\tree$ has a marked ray $\ray$, then for $v\in\ray$, $\Znorm_n^v=\ip{\vec Z^v_n}{\rv}/\rho^n$ is not necessarily a martingale for the first $\abs{h(v)}$ steps. Nevertheless it is eventually a martingale so we can still define $W_v$ to be the a.s.\ limit of $\Znorm^v_n$.} For vertices $v\in\tree$ we define $S_v$ as in \cite[\S3]{\PZ}: if $\tree$ is a rooted tree, let
\beq \label{e:harm.coords.rooted}
S_v \equiv \sum_{\rt < u \le v} W_u.
\eeq
If $\tree$ has marked ray $\ray$, recalling \eqref{e:horocycle} we set
\beq \label{e:harm.coords.ray}
S_v \equiv S_{R_v} + S^\ray_v\quad\text{where}\quad
S_{R_v} \equiv -\sum_{u\in\ray, \rt \ge u>R_v} W_u,\quad
S^\ray_v \equiv \sum_{R_v<u\le v} W_u.
\eeq
While on $\mgw$-a.e.\ $\tree$ the map $v\mapsto S_v$ is harmonic except at $\rt$ with respect to the transition probabilities of $\rw_\rho(\tree)$, on $\imgw$-a.e.\ $(\tree,\ray)$ the map $v\mapsto S_v$ is harmonic at every vertex with respect to the transition probabilities of $\rw_\rho(\tree,\ray)$. Thus, if $(Y_t)_{t\ge0}\sim\rw_\rho(\tree,\ray)$, $M_t\equiv S_{Y_t}$ will be a martingale given a fixed realization of the tree; we regard it as providing ``harmonic coordinates'' for the random walk. Using the reversing measure $\imgwr$ it is easy to prove a quenched CLT for $M$ (extending \cite[Cor.~1]{\PZ}):

\bppn
\label{p:mg.clt}
Under \hypirred, \hypks, and \hypmmt{2}, on $\imgw$-a.e.\ $(\tree,\ray)$ the process $M_{\flr{nt}}/(\eta\si\sqrt{n})$ converges in distribution to a standard Brownian motion as $n\to\infty$.

\bpf
We check the conditions of the Lindeberg--Feller martingale CLT (see e.g.\ \cite[Thm.~7.7.4]{durrett}): letting
$$V_n = \f{1}{n} \sum_{t=0}^{n-1} \E_\tree [(M_{t+1}-M_t)^2\giv \cG^\tree_t],$$
we verify that for $\imgw$-a.e.\ $(\tree,\ray)$,
\bnm[(i)]
\item $V_n\to \eta^2\si^2$ in probability and
\item for all $\ep>0$, $\f{1}{n} \sum_{t=0}^{n-1} \E_\tree[ (M_{t+1}-M_t)^2 \Ind{\abs{M_{t+1}-M_t} > \ep\sqrt{n}} ]\to0$.
\enm
Let $Y_n$ denote the random walk on $(\tree,\ray)$: we rewrite $V_n$ in terms of the induced random walk on $\treesray$ as
$$V_n = \f{1}{n}\sum_{t=1}^{n-1} \vph[(\tree,\ray)^{Y_t}],\quad
\vph[(\tree,\ray)]\equiv \f{\rho}{\rho+d_\rt} W_\rt^2 
	+ \f{1}{\rho+d_\rt}\sum_{j=1}^{d_\rt} W_{0j}^2.$$
By Thm.~\ref{t:imgw}~(\ref{t:imgw.erg}) and the Birkhoff ergodic theorem, we have $V_n$ converging $\imgwr$-a.s.\ to $\E_\imgwr[\vph]$ provided $\vph\in L^1(\imgwr)$. We calculate
$$\E_\imgwr[\vph]
= \f{1}{2\rho} \E_\mgw\bigg[
	\rho W_\rt^2 + \sum_{v\in\pd\rt} W_v^2
	\bigg]
= \E_\mgw[W_\rt^2]
= \eta^2\si^2,$$
so condition (i) is proved. Condition (ii) is checked similarly using dominated convergence.
\epf
\eppn

\begin{rmk}
To give some indication of how our results might be extended to $\rwre_{\bar\rho}$, we note that the main ingredient needed is the appropriate generalization of the normalized population size: we define it to be the random variable $\ol W_\rt$ which is the a.s.\ limit of the martingale $\Znorm_n\equiv \Znorm_n^{(1)}$ defined by \eqref{e:rwre.mg}. 
If $\ol W_v$ denotes the normalized population size of $\subtree{\tree}{v}$, then
$$\bar\rho \ol W_v = \sum_{w\in\desc v}\al_w \ol W_w,$$
so the $\ol W_v$ can be used to define harmonic coordinates for the RWRE. In the single-type case, $\ol W_\rt$ has finite second moment if and only if $\ka\ge2$ \cite[Thm.~2.1]{\Liu}, so clearly Propn.~\ref{p:mg.clt} cannot apply outside this regime. We emphasize again that due to the same technical barriers which arise in \cite{\Faraud}, simple adaptations of our proof will not cover the full regime $\ka\ge2$.
\end{rmk}

\subsection{Quenched $\imgwr$-CLT}
\label{ss:qiclt}

We now prove the quenched CLT for $\imgwr$ trees by controlling the corrector
$$\vep_t \equiv \f{M_t}{\eta}-h(Y_t)$$
on the interval $0\le t\le n$. For $1/2<\de<1$ and $n\ge 0$ fixed, let $\rtime{n}{j}$, for $j\flr{n^\de}\le n$ denote integer times chosen uniformly at random (independently of one another and of the random walk $Y$) from the interval $[j\flr{n^\de},(j+1) \flr{n^\de})$.

\bppn \label{p:corr.imgw}
Assume \hypirred, \hypks, and \hypmmt{\mmt} with $\mmt>2$. There exists
$\de_0\equiv\de_0(\mmt)\in(1/2,1)$ such that for $\de_0\le\de<1$ and $\ep>0$,
\beq \label{e:corr.imgw}
\lim_{n\to\infty} \P_{(\tree,\ray)}
	\Big(
	\max_{j\flr{n^\de}\le n} \absb{\vep_{\rtime{n}{j}}}
	\ge\ep\sqrt{n} \Big)
=0,\quad\imgwr\text{-a.s.}
\eeq
Further, for any $\ep'$ with $2\ep'+\de<1$,
\beq \label{e:tight.imgw}
\lim_{n\to\infty} \P_{(\tree,\ray)}
	\Big(
	\max_{r,s\le n, \abs{r-s}\le n^\de}
	\abs{h(Y_r)-h(Y_s)} \ge n^{1/2-\ep'}
	\Big)
=0,\quad\imgwr\text{-a.s.}
\eeq
\eppn

Given this proposition, we can prove the quenched CLT for $\rw_\rho$ on $\imgwr$ trees:

\bpf[Proof of Thm.~\ref{t:iclt}]
If $t\le n$ then $\abs{t-\rtime{n}{j}}\le \flr{n^\de}$ for some $j$, so
$$\max_{t\le n}\abs{\vep_t}
\le \max_{r,s\le n, \abs{r-s}\le \flr{n^\de}}
	\absb{ \f{M_r}{\eta}-\f{M_s}{\eta} }
+ \max_{j\flr{n^\de}\le n} \abs{ \vep_{\rtime{n}{j}} }
+ \max_{r,s\le n, \abs{r-s}\le \flr{n^\de}}
	\abs{ h(Y_r)-h(Y_s)}.$$
$M$ satisfies a CLT by Propn.~\ref{p:mg.clt}, and it follows from \eqref{e:corr.imgw} and \eqref{e:tight.imgw} that
$$\lim_{n\to\infty} \P_{(\tree,\ray)}
	\Big(
	\max_{t\le n} \abs{\vep_t} \ge \ep\sqrt{n}
	\Big) = 0,\quad\imgwr\text{-a.s.},$$
which gives the result.
\epf

The remainder of this section is devoted to the proof of Propn.~\ref{p:corr.imgw}.

\subsubsection{Tightness}

We begin by proving \eqref{e:tight.imgw}, using some {\it a priori} (annealed) estimates for $\rw_\rho$ coming from the Carne--Varopoulos bound.

\blem\label{l:carne}
There exists a constant $C<\infty$ such that
$$\mgw\Big(
	\max_{t\le n} \abs{X_t} \ge m \Big)
\le Cn e^{-(m+1)^2/(2n)}
\quad\forall m,n\ge1.$$

\bpf
We modify the proof of \cite[Lem.~5]{\PZ}. Take the finite tree with vertices $\set{w\in\tree : \abs{w} \le m}$, and make this into a wired tree $\tree^\star$ by adding a new vertex $\rt^\star$ which is joined by an edge to each vertex in $D_m$. Define the modified random walk $X^\star$ on $\tree^\star$ which follows the law of $\rw_\rho$ except at $\rt^\star$ where it moves to a vertex chosen uniformly at random from $D_m$. Then
$$\P_\tree(\max_{t\le n}\abs{X_t} \ge m) \\
\le 2 \sum_{t=1}^{n+1} \P_{\tree^\star}(X_t =\rt^\star).$$
By the Carne--Varopoulos inequality (see \cite[Thm.~13.4]{lyons.peres}),
$$\P_{\tree^\star}(X^\star_t = \rt^\star)
\le 2 \sqrt{\f{\abs{D_{m}}}{\rho^{m-1}}} e^{-(m+1)^2/(2t)}.
$$
Taking expectations gives
$$\mgw(\abs{X^\star_t} = \rt^\star)
\le C e^{-(m+1)^2/(2t)},$$
and summing over $1\le t\le n+1$ gives the result.
\epf
\elem

\bcor
\label{c:carne}
There exists a constant $C<\infty$ such that for any $m,n\ge1$,
$$C^{-1}\, \imgw_0\Big( \max_{t\le n} \abs{h(Y_t)}\ge m\Big)
\le \imgwr\Big( \max_{t\le n} \abs{h(Y_t)}\ge m\Big)
\le C n^2 e^{-m^2/(2n)}.$$

\bpf
We argue as in the proof of \cite[Cor.~2]{\PZ}. By decomposing into at most $n$ excursions away from height zero and using the stationarity of $\imgwr$, we find
\begin{align*}
&\imgwr \Big( \max_{t\le n} h(Y_t)\ge m \Big)  \\
&\le n \, \imgwr \Big( \exists t\le n :
	h(Y_t)\ge m, h(Y_s)>0 \ \forall 0\le s\le t \Big) \\
&\le C n \, \mgw\Big(
	\max_{t\le n} \abs{X_t} \ge m-1 \Big)
\le Cn^2 e^{-m^2/(2n)},
\end{align*}
by Lem.~\ref{l:carne}. The same bound holds for $\imgwr(\min_{t\le n} h(Y_t) \le -m)$ by the reversibility of $\imgwr$. The result follows by noting that $d\imgw_0/d\imgwr$ is uniformly bounded by a deterministic constant.
\epf
\ecor

\bpf[Proof of Propn.~\ref{p:corr.imgw}, \eqref{e:tight.imgw}]
By stationarity of $\imgwr$ and Cor.~\ref{c:carne}, for any fixed $s$
$$\imgwr \Big( \max_{0\le u\le n^\de} \abs{h(Y_{s+u})-h(Y_s)} \ge n^{1/2-\ep'}\Big)
\le C n^{2\de} e^{-n^{1-2\ep'-\de}/2},$$
and summing over $s\le n$ gives
$$\imgwr \Big( \max_{r,s\le t, \abs{r-s}\le n^\de}
	\abs{h(Y_r)-h(Y_s)} \ge n^{1/2-\ep'} \Big)
\le C n^{2\de+1} e^{-n^{1-2\ep'-\de}/2},$$
which is summable in $n$ provided $2\ep'+\de<1$. The result then follows from Markov's inequality and Borel--Cantelli.
\epf

\subsubsection{Control of corrector}
\label{sss:corrector}

In the remainder of this section we prove \eqref{e:corr.imgw}. We will make use of the following classical result:
\blem[{\cite[p.~60]{\Petrov}}]
\label{l:petrov}
If $z_1,\ldots,z_n$ are independent random variables with $\E z_i=0$ and $\E\abs{z_i}^\mmt<\infty$, then
$$\E\bigg[
	\bigg| \sum_{i=1}^n z_i \bigg|^\mmt\bigg]
\le \begin{cases}
2 \sum_{i=1}^n \E[\abs{z_i}^\mmt] & \text{if } 1\le\mmt\le2, \\
C(\mmt) n^{\mmt/2-1} \sum_{i=1}^n \E[\abs{z_i}^\mmt] & \text{if } \mmt\ge2.\end{cases}$$
\elem

Recalling \eqref{e:horocycle}~and~\eqref{e:harm.coords.ray}, we decompose
\beq \label{e:corr.imgw.decomp}
\f{1}{\sqrt{n}} \max_{j \flr{n^\de} \le n} \absb{\vep_{\rtime{n}{j}}}
\le E_1 + E_2
\eeq
where, with $R_t\equiv R_{Y_t}$ denoting the nearest ancestor of $Y_t$ on $\ray$,
$$E_1 \equiv \f{1}{\sqrt{n}} \max_{t\le 2n} \absb{\f{S_{R_t}}{\eta}-h(R_t)},\quad
E_2 \equiv \f{1}{\sqrt{n}} \max_{j \flr{n^\de} \le n}
		\absb{ \f{S^\ray_{Y_{\rtime{n}{j}}}}{\eta} 
			- d( Y_{\rtime{n}{j}},\ray)}.$$
The following lemma says that the harmonic coordinates $(S_v)_{v\in\tree}$ of \eqref{e:harm.coords.rooted}, rescaled by $\eta$ of \eqref{e:eta.sigma}, are a good approximation to the actual coordinates $\abs{v}$ on the $\mgw$ rooted trees. Let
\beq \label{e:mgw.coords.off.set}
A^\ep_n
\equiv A^\ep_n(\tree)
\equiv\setb{v\in D_n:\absb{\f{S_v}{n}-\eta}>\ep},\quad
\ep>0, \, n\ge1.
\eeq
Let $\Tret\equiv \min\set{t>0 : X_t=\rt}$ denote the first return time to the starting point $X_0=\rt$ by the walk $X$.

\blem\label{l:mgw.coords}
Assume \hypirred, \hypks, and \hypmmt{\mmt} with $\mmt\ge2$. For any $\ep>0$, the expected number of visits to $A^\ep_k$ during a single excursion away from the root is
$$\E_\mgw \bigg[ \sum_{t=0}^{\Tret} \Ind{X_t \in A^\ep_k} \bigg]
\le \f{C}{\rho^k} \E_\mgw \bigg[
	\sum_{v\in A^\ep_k} (1+d_v)
	\bigg]
\le \f{C(\mmt,\ep)}{k^{\mmt/2}}.$$

\bpf
If $v\in\tree$  with $\abs{v}=k\ge1$, a simple conductance calculation (see \cite[Ch.~2]{lyons.peres}) gives
\beq \label{e:mgw.exp.visits}
\E_\tree
	\bigg[ \sum_{t=0}^{\Tret} \Ind{X_t=v} \bigg]
= \f{\P_\tree(\rt\to v)}{\P_\tree(v\to\rt)}
= \f{\rho+d_v}{d_\rt \rho^k},
\eeq
so the first inequality follows. For the second we follow the proof of \cite[Lem.~3]{\PZ} (in particular the estimate \cite[(20)]{\PZ}) and of \cite[Lem.~4.2]{\Faraud}. Recall from \S\ref{ss:imgwr.wk.lim} the definition \eqref{e:Qn.a} of the probability measure $\Qntype{a}_k$ on rooted trees $\tree$ given by a size-biasing of $\mgwtype{a}$, and further the probability $\Qntype{a}_{k\star}$ on rooted trees $\tree$ with a marked path $(\rt=v_0,\ldots,v_k)$ from the root to level $k$:
\begin{align*}
&\E_{\mgwtype{a}}\bigg[
	\sum_{v\in A^\ep_k} (1+d_v) \bigg]
\le C \E_{\mgwtype{a}}\bigg[
	\sum_{v\in A^\ep_k} \rvec{\type_v} (1+d_v) \bigg] \\
&\le C\rho^k \E_{\Qntype{a}} \bigg[ \f{\sum_{v\in A^\ep_k} \rvec{\type_v} (1+d_v)}
	{\ip{\vec Z_k}{\rv}}\bigg]
= C \rho^k
	\E_{\Qntype{a}_{k\star}}
	[ (1+d_{v_k}) \Ind{v_k  \in A^\ep_k} ],
\end{align*}
so it suffices to show
$$\E_{\Qntype{a}_{k\star}}
	[(1+d_{v_k}) \Ind{v_k  \in A^\ep_k}]
\le C(\mmt,\ep) k^{-\mmt/2}.$$
To this end, writing $W_i\equiv W_{v_i}$, for $i<k$ we decompose $W_i \equiv W_{i+1}/\rho + W^\backprime_i$ where $W^\backprime_i$ is the normalized population size of $\subtree{\tree}{v_i}\setminus \subtree{\tree}{v_{i+1}}$. Then
$$W_i = \sum_{j=i}^{k-1} \f{W^\backprime_j}{\rho^{j-i}} + \f{W_k}{\rho^{k-i}},$$
so
$$\f{S_{v_k}}{k}-\eta = \f{1}{k} C_k W_k + \f{1}{k} \sum_{i=1}^{k-1} C_i W^\backprime_i -\eta,\quad C_i \equiv \sum_{j=0}^{i-1} \rho^{-j} \le C_\infty \equiv \f{\rho}{\rho-1}.$$
Conditional on the types $(\type_i \equiv \type_{v_i})_{i=1}^k$, the random variables $W^\backprime_1,\ldots,W^\backprime_{k-1}$ are independent of one another and of the pair $(W_k,d_{v_k})$, and all these random variables have finite moments of order $\mmt$ by Propn.~\ref{p:moments}. Therefore
$$\E_{\Qntype{a}_{k\star}}[(1+d_{v_k}) \Ind{v_k  \in A^\ep_k} ]
\le C \Qntype{a}_{k\star} \bigg(
	\bigg|
	\f{1}{k} \sum_{i=1}^{k-1} C_i W^\backprime_i -\eta
	\bigg|\ge \ep/2
	\bigg)
+ \E_{\Qntype{a}_{k\star}}
	[ (1+d_{v_k}) \Ind{ C_k W_k/k \ge\ep/2 } ]$$
By \hypmmt{\mmt}, Markov's inequality, and H\"older's inequality, the second term is
$$\le \bigg(
	\f{2 C_k}{k \ep} \bigg)^{\mmt-1}
	\E_{\Qntype{a}_{k\star}}[(1+d_{v_k}) W_k^{\mmt-1}]
\le \f{C(\mmt,\ep)}{k^{\mmt-1}} \le \f{C(\mmt,\ep)}{k^{\mmt/2}}$$
(since $\mmt\ge2$). As for the first term, by Lem.~\ref{l:petrov} and Markov's inequality,
\begin{align*}
&\Qntype{a}_{k\star}
\bigg( \bigg|
	\f{1}{k} \sum_{i=1}^{k-1} C_i W^\backprime_i
	-
	\E_{\Qntype{a}_{k\star}}\bigg[
	\f{1}{k} \sum_{i=1}^{k-1} C_i W^\backprime_i
	\,\bigg|\,
	(\type_i)_{i=1}^k
	\bigg]
	\bigg|
	>\ep/4 \bigg) \\
&\le C(\mmt) \f{k^{\mmt/2-1}}{(k\ep)^{\mmt}}
	\bigg\{ \sum_{i=1}^{k-1}
	\E_{\Qntype{a}_{k\star}}
	[\abs{ C_i (W^\backprime_i - \E[W^\backprime_i\giv \type_i])}^\mmt]
	 \bigg\}
\le C(\mmt,\ep) k^{-\mmt/2}.
\end{align*}
On the other hand,
$$\E_{\Qntype{a}_{k\star}}\bigg[
	\f{1}{k} \sum_{i=1}^{k-1} C_i W^\backprime_i
	\bigg]
\to C_\infty \E_{\Qn_\infty}[W^\backprime_\rt]
= \E_{\Qn_\infty}[W_\rt] = \eta,$$
and so
$$\Qntype{a}_{k\star}
\bigg( \bigg|
	\E_{\Qntype{a}_{k\star}}\bigg[
		\f{1}{k} \sum_{i=1}^{k-1} C_i W^\backprime_i
		\,\bigg|\,
		(\type_i)_{i=1}^k
	\bigg]
	-\eta
	\bigg|
	>\ep/4
\bigg)$$
decays exponentially in $k$ by \cite[Thm.~3.1.2]{\DZ}. Combining these estimates completes the proof.
\epf
\elem

Recalling the definition \eqref{e:harm.coords.ray} of the harmonic coordinates on the $\imgwr$ trees, the next step is to use Lem.~\ref{l:mgw.coords} to show that on these trees $S^\ray_v/\eta$ is a good approximation to $d(v,\ray)$. In analogy with \eqref{e:mgw.coords.off.set} set
\beq \label{e:imgw.coords.off.set}
B^\ep_k
= \setb{w\in \tree:d(w,\ray)=k, \absb{\f{S_w^\ray}{k}-\eta}>\ep},\quad
B^\ep = \bigcup_{k\ge1} B^\ep_k(\tree,\ray).
\eeq

\blem\label{l:imgw.coords}
Assume \hypirred, \hypks, and \hypmmt{\mmt} with $\mmt>2$. There exists $\de_0\equiv\de_0(\mmt)\in(1/2,1)$ such that for $\de_0\le\de<1$ and $\ep>0$,
$$\lim_{n\to\infty} \P_{(\tree,\ray)}
	\Big(
	\exists j\in\Z_{\ge0}, j\flr{n^\de}\le n : Y_{\rtime{n}{j}}\in B^\ep
	\Big)=0,\quad \imgwr\text{-a.s.}$$

\bpf
We modify the proof of \cite[(22)]{\Faraud}. If we define
$$\Thit_{n,\ep}
\equiv\inf\set{t\ge0:
	\abs{h(Y_t)}=\flr{n^{1/2+\ep}}},$$
then Cor.~\ref{c:carne} together with Markov's inequality gives
$$\imgwr[\P_{(\tree,\ray)}(\Thit_{n,\ep}\le n)\ge c]
\le c^{-1}\imgwr(\Thit_{n,\ep}\le n)
\le c^{-1} Cn^2 e^{-n^{2\ep}/2},$$
so by Borel--Cantelli we have $\P_{(\tree,\ray)}(\Thit_{n,\ep}\le n)\to0$, $\imgwr$-a.s.

On the event $\set{\Thit_{n,\ep}>n}$, we decompose the walk into excursions from $\ray$ started at $v_i$, $0\le i<\flr{n^{1/2+\ep}}$ (with each step of the walk along the ray contributing an empty excursion) and apply Wald's identity (see e.g.\ \cite[Exercise 22.8]{\Bill}) to find
\begin{align}
\nonumber &\P_{(\tree,\ray)}\Big(
	\Big\{
		\exists j\flr{n^\de}\le n :
		Y_{\rtime{n}{j}} \in B^\ep
		\Big\}
	\cap \set{\Thit_{n,\ep}>n}
	\Big)\\
\label{e:imgw.decomp.map}
&\le \f{1}{\flr{n^\de}} \sum_{i=0}^{\flr{n^{1/2+\ep}}-1}
	\E_{(\tree,\ray)}[
		\locT{\Thit_{n,\ep}}{i}
		]
	\,\E_{(\tree,\ray)}^i[
		\visits(B^\ep;\Texc)
	].
\end{align}
In the above, $\locT{n}{A}\equiv\visits(A;n)$ denotes the number of visits to set $A$ by time $n$ and $\locT{n}{i}\equiv\visits(v_i;n)$. $\E_{(\tree,\ray)}^i$ denotes expectation with respect to the law of a $\rho$-biased random walk $Y$ started from $Y_0=v_i$, and 
$\Texc\equiv\inf\set{t>0: Y_t=v_i \text{ or }
	Y_t\notin\subtree{\tree}{v_i}}$
denotes the excursion end time.

By a conductance calculation,
\beq\label{e:conduct.ray}
\E_{(\tree,\ray)}[\locT{\Thit_{n,\ep}}{i}]
=\f{1}{\P(v_i\to v_{\flr{n^{1+\ep}}})}
\le \f{1+d_i/\rho}{1-\rho} \le C d_i.
\eeq
During a single excursion away from $\ray$ the walk can visit only one of the $\subtree{\tree}{w}$ for $w\in\desc v \setminus v_{i-1}$, so to bound the second factor of each summand in \eqref{e:imgw.decomp.map} it suffices to consider an $\mgw$ rooted tree $\tree'$ (without ray): letting
$$\wt A^\ep_k
\equiv \wt A^\ep_k(\tree') \equiv \setb{v\in D_k : \absb{\f{W_\rt+S_v}{k+1}-\eta}>\ep},\quad \wt A^\ep \equiv \bigcup_{k\ge0} \wt A^\ep_k,$$
it follows from a (very slight) modification of Lem.~\ref{l:mgw.coords} that
$$\E_{\imgw_0}[
	\visits(B^\ep;\Texc)
	\giv i \cup \desc i] \le C
	  \E_\mgw \bigg[ \sum_{t=0}^{\Tret} \Ind{X_t \in \wt A^\ep} \bigg]
\le C(\mmt,\ep) \sum_{k\ge1} k^{-\mmt/2}
\le C(\mmt,\ep)$$
(using $\mmt>2$). It follows that the quantity in \eqref{e:imgw.decomp.map} converges to zero $\imgwr$-a.s., which concludes the proof.
\epf
\elem

\bpf[Proof of Propn.~\ref{p:corr.imgw}, \eqref{e:corr.imgw}]
Recall the decomposition \eqref{e:corr.imgw.decomp}. For any $k_0$,
$$ E_1 \le \f{1}{\sqrt{n}}\max_{i \le k_0} \absb{ \f{S_{v_i}}{\eta}-h(v_i) }
+ \f{1}{\sqrt{n}}\max_{t\le 2n, h(R_t) > k_0} \absb{\f{S_{R_t}}{\eta}-h(R_t)}.$$
The first term clearly tends to zero as $n\to\infty$ with $k_0$ fixed. The second term is bounded above by
\beq\label{e:dev.ray}
\bigg(
	\f{1}{\sqrt{n}} \max_{t\le 2n} \abs{M_t} \bigg)
	\sup_{i>k_0} \absb{\f{1}{\eta}-\f{h(v_i)}{S_{v_i}}}.
\eeq
Now recall from the proof of Propn.~\ref{p:imgw.wk} that if $(\tree,\ray)\sim\imgw_0$ then $(\subtree{\tree}{k},\rt)\sim\Qn_{k\star}$. Thus a consequence of the proof of Lem.~\ref{l:mgw.coords} is that for sufficiently small $\ep$,
$$\imgw_0(\abs{S_{v_k}/k+\eta}\ge\ep)
\le C(\mmt,\ep) k^{-\mmt/2}.$$
Therefore the supremum in \eqref{e:dev.ray} can be made arbitrarily small by taking $k_0$ large. We also have
$$E_2
\le \bigg(
	\f{1}{\sqrt{n}} \max_{t\le 2n} \abs{M_t} \bigg)
	\max_{j\flr{n^\de}\le n} \absb{\f{1}{\eta} - \f{d(Y_{\rtime{n}{j}},\ray)}{S^\ray_{Y_{\rtime{n}{j}}}} },$$
and in view of Lem.~\ref{l:imgw.coords} the second factor tends to zero in probability. By the invariance principle for $M$ proved in Propn.~\ref{p:mg.clt}, $\max_{t\le 2n} \abs{M_t}/\sqrt{n}$ stays bounded in probability as $n\to\infty$, so the result follows.
\epf

\section{\texorpdfstring{From $\imgwr$-CLT to $\mgw$-CLT by shifted coupling}
	{From IMGWR-CLT to MGW-CLT by shifted coupling}}
\label{s:coup}

In this section we prove our main result Thm.~\ref{t:mclt}. In \S\ref{ss:coup} we review (a slight modification of) the ``shifted coupling'' procedure of \cite[\S6]{\PZ}, which we use in \S\ref{ss:coup.amclt} to transfer the $\imgwr$-CLT to an annealed $\mgw$-CLT. In \S\ref{ss:coup.qmclt} we prove a variance estimate which allows to go from the annealed to the quenched $\mgw$-CLT.

\subsection{The shifted coupling construction}
\label{ss:coup}

We begin by reviewing the shifted coupling construction of \cite[\S6]{\PZ}, with the (natural) modification needed to handle the multi-type case. The basic observation underlying the construction is that the law of the random walk $X\sim\rw_\rho(\tree)$ up to time $t$ depends only on
$$\cE_t\equiv \rt \cup (\pd X_s)_{0\le s<t}$$
(``the subtree explored by time $t$''), so that one can construct the tree at the same time as the random walk.

For any tree $\tree$ (with or without marked ray) and $U$ any subset of the vertices of $\tree$, we also use $U$ to indicated the subgraph of $\tree$ induced by $U$. Let $\leaves\tree$ denote the set of leaves and $\tree^\circ\equiv\tree\setminus\leaves\tree$.

Let $a_0\in\cQ$ be fixed, and suppose $(\tree,(X_t)_{t\ge0})\sim\mgw\otimes\rw_\rho$. For each fixed $n\ge1$ we give a decomposition of $X$ into ``fresh excursions'' marked by time intervals $[\tau_i,\eta_i)$, $i\ge1$, as follows. Set $\eta_0\equiv0$ and define
\beq\label{e:ell}
\ell(n)\equiv4\flr{(\log (1+n))^{3/2}}.
\eeq
For $i\ge1$, let
\begin{align*}
\tau_i&\equiv\min\set{t>\eta_{i-1}:
	X_t\in\leaves\cE_t,
	\abs{X_t}>\ell(n)/2,
	\type_{X_t}=a_0
	}, &
	\text{excursion start,} \\
\eta_i&\equiv\min\set{t>\tau_i:
	X_t\in\cE_{\tau_i}^\circ}, &
	\text{excursion end,} \\
\cV_i&\equiv X_{\tau_i}\cup\cE_{\eta_i}\setminus\cE_{\tau_i}, &
	\text{excursion exploration.}
\end{align*}%
We take the convention $\min\emptyset\equiv\infty$, and let $\nxtime \equiv \max\setb{i:\eta_i<\infty}$ be the total number of excursions (so $\setb{\nxtime<\infty}=\nx$).

Next we construct a \emph{coupled} realization $((\tree^\coup,\ray),(Y_t)_{t\ge0})\sim\imgw_0\otimes\rw_\rho$ as follows: first construct the backbone $\cE^\coup_0$ of the tree ($\ray$ and $\desc v_i$ for $i\ge1$, together with types) in the manner described in \S\ref{ss:imgwr}. Set $\eta^\coup_0\equiv0$, and start a $\rho$-biased random walk $Y$ on $\cE^\coup_0$ with $Y_0=\rt$. As in the $\mgw$ setting we will construct a growing sequence $(\cE^\coup_t)_{t\ge0}$ such that $\cE^\coup_t=\cE^\coup_0 \cup (\pd Y_s)_{0\le s<t}$, and we will define (for $i\ge1$)
\begin{align*}
\tau^\coup_i&\equiv\min\set{t>\eta^\coup_{i-1}:
	Y_t\in\leaves\cE^\coup_t,
	d(Y_t,\ray)>\ell(n)/2,
	\type_{Y_t}=a_0}, &
	\text{excursion start,} \\
\eta^\coup_i&\equiv\min\set{t>\tau^\coup_i:
	Y_t\in(\cE^\coup_{\tau^\coup_i})^\circ}, &
	\text{excursion end,} \\
\cV^\coup_i&\equiv
	Y_{\tau^\coup_i}
	\cup\cE^\coup_{\eta^\coup_i} \setminus \cE^\coup_{\tau^\coup_i}, &
	\text{excursion exploration.}
\end{align*}
The difference is that we grow the sequence $\cE^\coup_t$ in a manner \emph{dependent} on $(\tree,X)$, such that excursions of $Y$ into unexplored territory (and started from $a_0$) match the excursions of $X$ defined above: formally, we couple $(Y_s)_{\tau^\coup_i\le s<\eta^\coup_i}$ with $(X_s)_{\tau_i\le s<\eta_i}$ such that there is a (type-preserving) isomorphism $f_i:\cV_i\to \cV^\coup_i$ with $f_i(X_{\tau_i+s})=Y_{\tau^\coup_i+s}$, and then we set $Y_{\eta^\coup_i}$ to be the ancestor of $Y_{\tau^\coup_i}$ (not necessarily of the same type as $X_{\eta_i}$). Then, on the inter-excursion intervals $\eta^\coup_{i-1}\le t < \tau^\coup_i$ (for $i\ge1$),
\bitm
\item If $Y_t\in(\cE^\coup_t)^\circ$ then generate $Y_{t+1}$ according to the transition kernel of $\rw_\rho$ on $\cE^\coup_{t+1}=\cE^\coup_t$;
\item If $Y_t\in\leaves\cE^\coup_t$ with $\type_{Y_t}\ne a_0$, let $\cE^\coup_{t+1}$ be the enlargement of $\cE^\coup_t$ obtained by attaching random offspring to $Y_t$ according to law $\q^{\type_{Y_t}}$, and generate $Y_{t+1}$ according to the transition kernel of $\rw_\rho$ on $\cE^\coup_{t+1}$.
\eitm
Finally, with $\cE^\coup_\infty\equiv\lim_{t\to\infty}\cE^\coup_t$, we define $\tree^\coup$ by attaching to each vertex $v\in\leaves\cE^\coup_\infty$ an independent $\mgw^{\type_v}$ tree. We thus obtain the following extension of \cite[Lem.~8]{\PZ}:

\blem\label{l:shifted.coupling}
If $(\tree,(X_t)_{t\ge0})\sim\mgw\otimes\rw_\rho$ then the marginal law of $((\tree^\coup,\ray), (Y_t)_{t\ge0})$ arising from the above construction is $\imgw_0\otimes\rw_\rho$.
\elem

\brmk\label{r:param.n}
Although we suppress the parameter $n$ from the notation, we emphasize that each $n\ge1$ gives rise to a different excursion decomposition, hence a different coupling between $(\tree,X)$ and $((\tree^\coup,\ray),Y)$.
\ermk

\subsection{Annealed \texorpdfstring{$\mgw$-CLT}{MGW-CLT}}
\label{ss:coup.amclt}

We now transfer the quenched $\imgwr$-CLT to the following annealed $\mgw$-CLT:

\bppn\label{p:mgw.aclt}
Assume \hypirred, \hypks\ and \hypmmt{\mmt} with $\mmt>4$. If $(\tree,X)$ has law $\mgw\otimes\rw_\rho$ conditioned on $\nx^c$, then the processes $(\abs{X_{\flr{nt}}}/(\si\sqrt{n}))_{t\ge0}$ converge in law to the absolute value of a standard Brownian motion.
\eppn

Recall that $R_t\equiv R_{Y_t}$ denotes the nearest ancestor of $Y_t$ on $\ray$. By Thm.~\ref{t:iclt}, for $\imgw_0$-a.e.\ $(\tree,\ray)$, the process
\beq\label{e:refl}
\reflec_{\flr{nt}}/(\si\sqrt{n}),\quad
\reflec_t\equiv
h(Y_t)-\min_{0\le s\le t} h(Y_s)
=h(Y_t)-\min_{0\le s\le t} h(R_s)
\ge0
\eeq
converges to a Brownian motion minus its running minimum, which is the same in law as the absolute value of a Brownian motion (see e.g.\ \cite[Thm.~3.6.17]{\KS}). Thus to deduce Propn.~\ref{p:mgw.aclt} we need to estimate the relation between the processes $\abs{X_n}$ and $\reflec_n$. To this end, let $\bt$, $\bt^\coup$ be the monotone increasing bijections
$$\bt:\Z_{\ge0} \to \bigcup_{i\ge1} [\eta_{i-1},\tau_i),\quad
\bt^\coup:\Z_{\ge0} \to\bigcup_{i\ge1} [\eta^\coup_{i-1},\tau^\coup_i),$$
parametrizing the inter-excursion times of $X_n$ and $Y_n$ respectively. We make the following notations (the left column refers to the $\mgw$ tree, while the right column refers to the $\imgw_0$ tree):
$$\begin{array}{ll}
X^\inter_s \equiv X_{\bt(s)},
& Y^\inter_s \equiv Y_{\bt^\coup(s)} \\
\cH_s \equiv \si(X_t : t\le\bt(s) )
& \cH^\coup_s \equiv \si(Y_t : t\le\bt^\coup(s) ), \\
J_i \equiv \bt^{-1}[\eta_{i-1},\tau_i),
& J^\coup_i \equiv (\bt^\coup)^{-1}[\eta^\coup_{i-1},\tau^\coup_i); \\
I_n\equiv\max\setb{i:\eta_{i-1}\le n},
& I^\coup_n\equiv\max\setb{i:\eta^\coup_{i-1}\le n}; \\
\De_n \equiv \sum_{i=1}^{I_n} \abs{J_i},
& \De^\coup_n \equiv \sum_{i=1}^{I^\coup_n} \abs{J^\coup_i}; \\
\De_n(\al) \equiv
	\sum_{i=1}^{I_n} \abs{\set{s \in J_i
	: \abs{X^\inter_s}\le n^\al}},
& \De^\coup_n(\al) \equiv \sum_{i=1}^{I^\coup_n}
	\abs{\set{s \in J^\coup_i
	: d(Y^\inter_s,\ray)\le n^\al}}.
\end{array}$$
In words, given the walk $X$ on the $\mgw$ tree, $X^\inter_s$ is the ``inter-excursion process'' adapted to the filtration $\cH_s$, $J_i$ is the $i$-th inter-excursion interval, $I_n$ is the number of such intervals intersecting $[0,n]$, $\De_n$ is the total length of these intervals, and $\De_n(\al)$ is the length of these intervals except for times spent at distance more than $n^\al$ from the root. The right column defines the analogous objects for the walk on the $\imgw_0$ tree.

\blem\label{l:azuma}
Assume \hypirred, \hypks, and \hypmmt{\mmt} with $\mmt>4$. There exists $\al_0(\mmt)<1/2$ such that for $\al>\al_0(\mmt)$,
$$\left.\begin{array}{r}
	\mgw(\De_n(\al)\ne\De_n \giv \nx^c) \\
	\imgw_0(\De^\coup_n(\al)\ne\De^\coup_n)
	\end{array}\right\} \le n^{-c},\quad
	c\equiv c(\mmt,\al)>0.$$
\elem

We will obtain the corollary below as a relatively straightforward consequence of Lem.~\ref{l:azuma}. 
Let
$$\drift_n
\equiv \max_{0\le r\le s\le n}
\set{h(R_s) - h(R_r)}$$
denote the maximum displacement by time $n$ against the backward drift on $\ray$.

\bcor\label{c:azuma}
Assume \hypirred, \hypks, and \hypmmt{\mmt} with $\mmt>4$. Then
\bnm[(a)]
\item\label{c:azuma.int} There exists $\al_0(\mmt)<1/2$ such that for $\al\ge\al_0(\mmt)$,
$$\left.
\begin{array}{r}
\mgw(\De_n \ge n^{1/2+\al+\ep} \giv \nx^c) \\
\imgw_0(\De^\coup_n \le n^{1/2+\al+\ep})
\end{array}
\right\}\le n^{-c},
\quad c\equiv c(\mmt,\al,\ep)>0$$
\item\label{c:azuma.drift} On $\imgw_0$-a.e.\ $(\tree,\ray)$, $\drift_n/\sqrt{n}$ converges $\P_{(\tree,\ray)}$-a.s.\ to zero.
\enm
\ecor

Assuming these results we can prove the annealed $\mgw$-CLT:

\bpf[Proof of Propn.~\ref{p:mgw.aclt}]
Let $\bb:\Z_{\ge0}\to\Z_{\ge0}$ be any nondecreasing map which maps $[\tau^\coup_i,\eta^\coup_i)$ bijectively onto $[\tau_i,\eta_i)$, and $J^\coup_i$ into $J_i$, for each $i$. Then for $t\in[\tau^\coup_i,\eta^\coup_i)$ we have
$\abs{X_{\bb(t)}}-\abs{X_{\tau_i}} = d(Y_t,\ray)-d(Y_{\tau^\coup_i},\ray)$, so, recalling \eqref{e:horocycle} and \eqref{e:refl}, we have
$$\abs{\abs{X_{\bb(t)}}-\reflec_t}
=\absb{
\abs{X_{\bb(t)}}-d(Y_t,\ray)-h(R_t)
	+\min_{s\le t} h(R_s)
}
\le\abs{X_{\tau_i}}+d(Y_{\tau^\coup_i},\ray)+\drift_t.
$$
If instead $t\in J^\coup_i$ then
$$\abs{\abs{X_{\bb(t)}}-\reflec_t}
\le \abs{X_{\bb(t)}}+\abs{\reflec_t}
\le \abs{X_{\bb(t)}}+\drift_t+d(Y_t,\ray).$$
It follows that on the event $\set{\De_n(\al)=\De_n}\cap\set{\De^\coup_n(\al)=\De^\coup_n}$,
$$\f{1}{\sqrt{n}} \max_{0\le t\le n}\abs{\abs{X_{\bb(t)}}-\reflec_t} \le
	\f{2n^\al + \drift_n}{\sqrt{n}},$$
so by Thm.~\ref{t:iclt}, Lem.~\ref{l:shifted.coupling}, Lem.~\ref{l:azuma}, and Cor.~\ref{c:azuma}~(\ref{c:azuma.drift}), the processes $(\abs{X_{\bb(\flr{nt})}}/(\si\sqrt{n}))_{t\ge0}$ converge in law to a reflected Brownian motion. On the other hand, Cor.~\ref{c:azuma}~(\ref{c:azuma.int}) implies that $n^{-1}\max_{0\le t\le 1} (\bb(\flr{nt})-\flr{nt})\to0$ in probability, so we obtain the CLT for the processes $(\abs{X_{\flr{nt}}}/(\si\sqrt{n}))_{t\ge0}$ from the a.s.\ uniform continuity of Brownian motion on compact intervals.
\epf

In the remainder of this subsection we prove Lem.~\ref{l:azuma} and Cor.~\ref{c:azuma}. Let $\cond{\rt}{\ell}\equiv \cC(\rt\leftrightarrow\ell)$ denote the conductance between $\rt$ and $D_\ell$ in $\tree$, with respect to the stationary measure $\vpi$ for $\rw_\rho(\tree)$ with the normalization $\vpi(\rt)=d_\rt$. We will make use of the following conductance lower bound:

\blem\label{l:polar}
Under \hypirred, \hypks, and \hypmmt{2}, there exist $0<r,C<\infty$ such that for all $\ep>0$,
$\mgw(\cond{\rt}{k}^{-1} \ge k^{1+\ep}\giv \nx^c) \le C k^{-r\ep}$.
\elem

The proof of the lemma is deferred to \S\ref{ss:appx.cond} where we also provide a quenched conductance lower bound (see Propn.~\ref{p:polar}) which is not needed in the proof of the main theorem. Lem.~\ref{l:polar} readily implies an upper bound on the amount of time
$$\Troot{n}{\al} \equiv \sum_{t=0}^n \Ind{\abs{X_t} \le n^\al},\quad
\Trootcoup{n}{\al}
	\equiv\sum_{t=0}^n \Ind{d(Y_t,\ray)\le n^\al}
$$
spent by $X$ (resp.\ $Y$) within distance $n^\al$ of the root (resp.\ marked ray) by time $n$:

\bcor\label{c:timeroot}
Assume \hypirred, \hypks, and \hypmmt{2}. Then
$$
\left.\begin{array}{r}
\mgw(\Troot{n}{\al} \ge n^{1/2+\al+\ep}\giv \nx^c) \\
\imgw_0(\Trootcoup{n}{\al} \ge n^{1/2+\al+\ep})
\end{array}\right\} \le Cn^{-c\ep}.
$$

\bpf
By iterated expectations, Markov's inequality, and Lem.~\ref{l:carne},
\begin{align*}
&\mgw(\Troot{n}{\al} \ge n^{1/2+\al+2\ep}\giv \nx^c)
\le n^{-\ep}
+\mgw(\E_\tree[\Troot{n}{\al}] \ge n^{1/2+\al+\ep}\giv \nx^c)\\
&\le n^{-\ep}
+ C e^{-n^{2\ep}/3}
+ \mgw( \set{\Thit_{n,\ep}>n} \cap 
	\set{
	\E_\tree[\Troot{\Thit_{n,\ep}}{\al}]
		\ge n^{1/2+\al+\ep}} \giv \nx^c)
\end{align*}
where $\Thit_{n,\ep}\equiv\inf\set{t\ge0: \abs{X_t}=\flr{n^{1/2+\ep}}}$. Wald's identity gives
$$\E_\tree[\Troot{\Thit_{n,\ep}}{\al}]
\le \E_\tree[\locT{\Thit_{n,\ep}}{\rt}]\,
\E_\tree[\Troot{\Tret}{\al}]$$
with $\locT{n}{A}$ the number of visits the walk makes to set $A$ by time $n$. Recalling \eqref{e:mgw.exp.visits}, we have
$$\E_\tree[\locT{\tau}{\rt}]
\le \f{d_\rt}
	{\cC(\rt\leftrightarrow\flr{n^{1/2+\ep}})},\quad
\E_\tree[N_{\Tret}(\al)]
\le C \sum_{k=0}^{\flr{n^\al}+1} \abs{D_k},$$
so the $\mgw$ bound follows from Lem.~\ref{l:polar} with a few more applications of Markov's inequality.

For the $\imgw_0$ bound we argue as in the proof of Lem.~\ref{l:imgw.coords}: by Markov's inequality and Cor.~\ref{c:carne},
\begin{align*}
&\imgw_0(\Trootcoup{n}{\al}
	\ge n^{1/2+\al+2\ep})
\le n^{-\ep}
	+ \imgw_0(
		\E_{(\tree,\ray)}
			[\Trootcoup{n}{\al}]
		\ge n^{1/2+\al+\ep} ) \\
&\le n^{-\ep} + C e^{-n^{4\ep}/3}
	+ \imgw_0(\Trootcoup{\Thit_{n,\ep}}{\al}
	\ge n^{1/2+\al+2\ep}),
\end{align*}
so it suffices to bound the last term. By Wald's identity and \eqref{e:conduct.ray},
\begin{align*}
\E_{(\tree,\ray)}
	[\Trootcoup{\Thit_{n,\ep}}{\al}]
&\le \sum_{i=0}^{\flr{n^{1/2+\ep}}-1}
	\E_{(\tree,\ray)}[\locT{\Thit_{n,\ep}}{i}]\,
	\E_{(\tree,\ray)}^i
		[ \Trootcoup{\Texc}{\al} ]\\
&\le C \sum_{i=0}^{\flr{n^{1/2+\ep}}-1} d_i \,
	\E_{(\tree,\ray)}^i[ \Trootcoup{\Texc}{\al} ],
\end{align*}
so again the bound follows by using Markov's inequality.
\epf
\ecor

Most of the technical estimates required for the proof of Lem.~\ref{l:azuma} are contained in the following auxiliary lemma (cf.\ \cite[Lem.~7.3]{\Faraud}). For $\ell(n)$ as in \eqref{e:ell}, define the sequence of $(\cH_s)$-stopping times
$$\Thet_0\equiv0,\quad
\Thet_{j+1} \equiv \min\set{s>\Thet_j : \abs{\abs{X^\inter_s}-\abs{X^\inter_{\Thet_j}} = \ell(n)}
}$$
and similarly the sequence of $(\cH^\coup_s)$-stopping times
$$\Thet^\coup_0\equiv0,\quad
\Thet^\coup_{j+1}
	\equiv\min
	\set{s>\Thet^\coup_j : \abs{d(Y^\inter_s,\ray)
	-d(Y^\inter_{\Thet^\coup_j},\ray)} = \ell(n)}.$$

\blem\label{l:coupling.est}
Assume \hypirred\ and \hypks.
\bnm[(a)]
\item\label{l:coupling.est.trunc}
Assume \hypmmt{\mmt} and let
$$\badC\equiv\badC(n,\ep)
\equiv\set{v\in\tree:W_v>n^{1/4-\ep}}$$
(well-defined for trees with and without ray). Then
$$\left.\begin{array}{r}
	\mgw(\Thit(\badC)\le n) \\
	\imgw_0(\Thit(\badC\setminus\ray)\le n)
	\end{array}\right\}
	\le C(\mmt,\ep) n^{1-\mmt(1/4-\ep)}.
$$
For $\mmt>4$ the right-hand side can be made $\le n^{-c}$ for $c\equiv c(\mmt,\ep)>0$.

\item\label{l:coupling.est.exc}
Assuming \hypmmt{\mmt} with $\mmt>2$, for any $\ep>0$ there exists $c\equiv c(\mmt,\ep)$ such that
$$\left.\begin{array}{r}
	\mgw(I_n \ge n^{1/2+\ep}\giv \nx^c) \\
	\imgw_0(I^\coup_n \ge n^{1/2+\ep} )
	\end{array}\right\}
	\le e^{-cn^{\ep/2}}.
$$

\item\label{l:coupling.est.samp}
With $\Thet_j,\Thet^\coup_j$ defined as above,
$$\left.\begin{array}{r}
\mgw(\bt(\Thet_{3I_n})\le n) \\
\imgw_0(\bt^\coup(\Thet^\coup_{3I^\coup_n})\le n)
\end{array}\right\}
\le e^{-c \ell(n)}.$$

\item\label{l:coupling.est.harm}
Recalling the notation of \eqref{e:mgw.coords.off.set} and \eqref{e:imgw.coords.off.set}, define
$$\badA\equiv\badA(n,\al,\ep)\equiv
	\bigcup_{k=0}^{\flr{(\log n)^2}}
	A^\ep_{\flr{n^\al}-k},\quad
\badB\equiv\badB(n,\al,\ep)\equiv
	\bigcup_{k=0}^{\flr{(\log n)^2}}
	B^\ep_{\flr{n^\al}-k}.$$
Assuming \hypmmt{\mmt} with $\mmt>2$, there exists $\al_0(\mmt)\in(0,1/2)$ such that for all $\al\ge\al_0(\mmt)$ and all $\ep<\ep_0(\mmt,\al)$ (with $\ep_0(\mmt,\al)>0$),
\beq\label{e:badset}
\left.\begin{array}{r}
\mgw(\Thit_\badA\le n) \\
\imgw_0(\Thit_\badB\le n)
\end{array}\right\}
\le n^{-c}
\eeq
for some $c\equiv c(\mmt,\al,\ep)>0$. If further $\mmt>4$ then $\al_0,\ep_0$ can be chosen such that
\beq \label{e:mgw.badset}
\mgw(\Thit_{\badA_+}\le n)
\le n^{-c},\quad
\badA_+\equiv\badA_+(n,\al,\ep)
\equiv\bigcup_{k\ge \flr{n^\al}} A^\ep_k.
\eeq
\enm

\bpf
(\ref{l:coupling.est.trunc}) See proof of \cite[(63)]{\PZ}.

\medskip\noindent
(\ref{l:coupling.est.exc}) We will show that with probability $\ge 1-e^{-cn^{\ep/2}}$ conditioned on $\nx^c$, one of the first $\flr{n^{1/2+\ep}/2}$ excursions has length $\eta_i-\tau_i>n$ which certainly implies the result. Conditioning on $\nx^c$ is needed simply to ensure $\nxtime=\infty$; for the purpose of proving the claim we may artificially define $\eta_i-\tau_i=\infty$ for $i> \nxtime$.

Then, conditioned on $(\eta_j-\tau_j)_{j=1}^{i-1}$, the probability that $\eta_i-\tau_i>n$ is bounded below by a constant times $\mgw(\Tret>n)$. Further
$$\P_\tree(\Tret>n)
\ge \P_\tree(\Tret>\Thit_{n,\ep/2}>n)
\ge \P_\tree(\Tret>\Thit_{n,\ep/2}) - \P_\tree(\Thit_{n,\ep/2}\le n),$$
so Lem.~\ref{l:carne} and Lem.~\ref{l:polar} imply
$\mgw(\Tret>n)\ge c/n^{1/2+\ep/2}$ for $c\equiv c(\mmt,\ep)$.
Thus the probability that none of the first $\flr{n^{1/2+\ep}}$ excursions has length $>n$ is
$$\le \Big(
	1-\f{c}{n^{1/2+\ep/2}}\Big)^{
		\flr{n^{1/2+\ep}/2}
	}
	\le e^{-c n^{\ep/2}},$$
which proves the result.

\medskip\noindent
(\ref{l:coupling.est.samp})
We follow the proof of \cite[Lem.~11]{\PZ}. On the $\imgw_0$ tree, since $d(Y_0,\ray)=0$, $d(Y^\inter,\ray)$ must increase by $\ell(n)$ going from $\Thet^\coup_{j-1}$ to $\Thet^\coup_j$ for at least half of the indices $j\le 3I^\coup_n$ so $\set{\bt(\Thet_{3I_n})\le n}$ implies the event
$$\badThet^\coup\equiv\set{
	\exists i\le I^\coup_n,
	\Thet^\coup_{j-1}, \Thet^\coup_j \in J^\coup_j,
	d(Y^\inter_{\Thet^\coup_j},\ray)
		> d(Y^\inter_{\Thet^\coup_{j-1}},\ray)
}.$$
This in turn implies one of two possibilities:
\bnm
\item[1.] there exist times $t_0<t_1<t_2\le n$ with $Y_{t_0}=Y_{t_2}$ and
$d(Y_{t_0},\ray) = d(Y_{t_1},\ray) + \ell(n)/4$, or
\item[2.] there exist times $t_1<t_2\le n$ with $d(Y_{t_2},\ray)=d(Y_{t_1},\ray)+\ell(n)/4$ such that $a_0$ does not appear on the geodesic between the $Y_{t_i}$.
\enm
By a random walk estimate (cf.\ \eqref{e:mgw.exp.visits}) summed over at most $n^2$ possibilities for $(Y_{t_0},Y_{t_1})$, the first event has probability $\le C n^2 \rho^{-\ell(n)/4}$. The second event has probability $\le e^{-c\ell(n)}$ by the construction of $\imgw_0$ and the irreducibility of the Markov chain, and combining these estimates gives the bound for $\imgw_0$. The bound for $\mgw$ follows by a similar argument.

\medskip\noindent
(\ref{l:coupling.est.harm})
We first prove the bounds for $\mgw$; the argument is similar to that of Cor.~\ref{c:timeroot}: again it suffices to bound $\mgw(\Thit_\badA\le\Thit_{n,\ep})$, and Wald's identity gives
$$\P_\tree(\Thit_\badA\le\Thit_{n,\ep})
\le \E_\tree[
	\locT{\Thit_{n,\ep}}{\rt}]\,
	\P_\tree[\Thit_\badA<\Tret].$$
By \eqref{e:mgw.exp.visits} and Markov's inequality, $\E_\tree[\locT{\Thit_{n,\ep}}{\rt}] \le n^{1/2+2\ep}$ except with probability at most $n^{-c}$ for $c\equiv c(\mmt,\ep)>0$. By Lem.~\ref{l:mgw.coords},
$$\mgw(\Thit_\badA<\Tret)
\le C(\al,\mmt,\ep) (\log n)^2 (n^\al)^{-\mmt/2}.$$
Since $\mmt>2$, we can choose $\al$ sufficiently close to $1/2$ and $\ep$ sufficiently small such that Markov's inequality gives $\mgw( \P_\tree(\Thit_\badA<\Tret) \ge n^{-(1/2+3\ep)})\le n^{-c}$, from which \eqref{e:badset} follows for $\mgw$. \eqref{e:mgw.badset} follows similarly by noting
$$\mgw ( \Thit_{\badA_+}
	<\Tret)
	\le C(\al,\mmt,\ep) (n^\al)^{1-\mmt/2}.$$
The bound \eqref{e:badset} for $\mgw$ together with the argument of Cor.~\ref{c:timeroot} gives
$$\imgw_0(\Thit_\badB\le\Thit_{n,\ep})\le C(\mmt,\ep) (\log n)^2 n^{1/2+\ep-\al\mmt/2},$$ and the bound \eqref{e:badset} for $\imgw_0$ follows by choosing $\al$ close to $1/2$ and $\ep$ small.
\epf
\elem

\bpf[Proof of Lem.~\ref{l:azuma}]
We modify the proof of \cite[(55)]{\PZ} (see also \cite[(28)]{\Faraud}; our Lem.~\ref{l:coupling.est} plays the role of \cite[Lem.~5.1]{\Faraud}).

On the $\mgw$ tree write $\badE\equiv\set{\De_n(\al)\ne\De_n}=\set{\max_{s\le\De_n} \abs{X^\inter_s}\ge n^\al}$. If we define
$$\badY\equiv\set{I_n\ge n^{1/2+\ep}}
\cup\set{\bt(\Thet_{3I_n}\le n)}
\cup\set{\Thit_\badA\le n}
\cup\set{\Thit_\badC\le n}$$
and consider the process $\cM(j)\equiv S_{X^\inter(\Thet_j)}$, we have
\beq\label{e:badE.decomp}
\mgw(\badE)
\le \mgw(\badY)	
	+ \mgw \Big(
		\badY^c \cap \Big\{
		\max_{j\le 3n^{1/2+\ep}}
		\cM(j) \ge
		(\eta-\ep) (n^\al-\ell(n))
		\Big\} \Big).
\eeq
By Lem.~\ref{l:coupling.est} it suffices to bound the second term. To this end let $\Tret_i$ denote the $i$-th return of $(X^\inter_{\Thet_j})_j$ to the root (with $\Tret_0\equiv0$): then for each $i\ge0$ the process
$\cM_i(j)\equiv
	\cM(j\wedge\Tret_{i+1})$
is a supermartingale for $j\ge\Tret_i+1$, and the second term of \eqref{e:badE.decomp} is
$$\le
\mgw \Big(
		\badY^c \cap \Big\{
		\max_{i:\Tret_i\le 3n^{1/2+\ep}}
		\max_{j\le 3n^{1/2+\ep}}
		[\cM_i(\Tret_i+j)-\cM_i(\Tret_i+1)] \ge
		\f{\eta n^\al}{2}
		\Big\} \Big)
$$
(for $n$ large and suitable $\al,\ep$). If we define the $(\cH_{\Thet_j})_j$-stopping time
$$\Psi\equiv\inf\set{j:\bt(\Thet_j)>\Thit_{\badC}},$$
then
$$\mathfrak{M}_i(j)\equiv
	\cM_i(j\wedge\Psi)
	-[\cM_i(\Psi)-\cM_i(\Psi-1)] \Ind{\Psi\le j},\quad
	j\ge\Tret_i+1$$
is a supermartingale with differences $\le\ell(n)n^{1/4-\ep}$, so the Azuma--Hoeffding inequality gives
$$\P_\tree\Big(
	\max_{j\le 3n^{1/2+\ep}} [\mathfrak{M}_i(\Tret_i+j)-\mathfrak{M}_i(\Tret_i+1)]
	\ge \f{\eta n^\al}{2}
	\Big)
\le \exp\ls
	-\f{ (\eta n^\al/2)^2 }
		{2(3n^{1/2+\ep}) \ell(n)^2 n^{2(1/4-\ep)}}
	\rs.$$
Choosing $\al,\ep$ appropriately and summing over at most $3n^{1/2+\ep}$ return times $\Tret_i$ gives the desired bound on the second term of \eqref{e:badE.decomp}, from which the $\mgw$ bound follows.

The bound for the $\imgw_0$-probability of $\badE^\coup\equiv\set{\De^\coup_n(\al)\ne\De^\coup_n}$ is similar, indeed simpler since $\cM^\coup(j)\equiv S_{Y^\inter(\Thet^\coup_j)}$ is always a supermartingale. With
$$\badY^\coup\equiv
	\set{I^\coup_n\ge n^{1/2+\ep}}
	\cup\set{\bt^\coup(\Thet^\coup_{3I^\coup_n})\le n}
	\cup\set{\Thit_\badB\le n}
	\cup\set{\Thit_\badC\le n},$$
we have from Lem.~\ref{l:coupling.est} that $\imgw_0(\badE^\coup)$ is
$$\le n^{-c}
	+\imgw_0\Big(
		(\badY^\coup)^c
		\cap\Big\{
			\max_{j\le 3n^{1/2+\ep}}
			\cM^\coup(j) \ge \f{\eta n^\al}{2}
			\Big\}
	\Big),$$
and applying the Azuma--Hoeffding bound gives the result.
\epf

\bpf[Proof of Cor.~\ref{c:azuma}]
\medskip\noindent
(\ref{c:azuma.int})
We have the set inclusions
\begin{align*}
\set{\De_n\ge n^{1/2+\al+\ep}}
&\subseteq \set{\De_n \ne\De_n(\al)}
	\cup\set{\Troot{n}{\al}\ge n^{1/2+\al+\ep}},\\
\set{\De^\coup_n\le n^{1/2+\al+\ep}}
&\subseteq\set{\De^\coup_n\ne\De^\coup_n(\al)}
	\cup\set{\Trootcoup{n}{\al} \ge n^{1/2+\al+\ep}},
\end{align*}
so the result follows from Lem.~\ref{l:azuma} and Cor.~\ref{c:timeroot}.

\medskip\noindent
(\ref{c:azuma.drift}) Let $(h_s)_{s\ge0}$ denote the height process for the walk $Y$ restricted to $\ray$, i.e.\ erasing all excursions away from $\ray$; clearly $\drift_n\le\drift'_n\equiv\max\set{h_s-h_r : 0\le r\le s\le n}$. But $h_s$ is simply a random walk on $\Z_{\le0}$ with a $\rho$-bias in the negative direction. Set $\si_0\equiv0$,
$$\si_j\equiv\inf\set{s>\si_{j-1}
	: h_s = h_{\si_{j-1}}-1}, \quad j\ge1.$$
Now the processes $(\bar h^{(j)}_s\equiv h_s-h_{\si_j})_{\si_j\le s\le\si_{j+1}}$ are i.i.d., and clearly $\si_n\ge n$, so
$$\drift_n'\le\max_{0\le j<n} \left(\max_s \bar h^{(j)}_s\right).$$
The probability of $\max_s \bar h^{(j)}_s\ge m$ is at most the probability that a random walk on $\Z$ started at $0$ with a $\rho$-bias in the negative direction will reach $m$ before $-1$, which is $(1-\rho^{-1})/(\rho^m-\rho^{-1}) \le \rho^{-m}$. Summing over $j$ gives $\P_{(\tree,\ray)}(\drift_n'\ge m)\le n\rho^{-m}$, $\imgw_0$-a.s.
\epf

\subsection{Quenched \texorpdfstring{$\mgw$-CLT}{MGW-CLT}}
\label{ss:coup.qmclt}

We now describe how to move from the annealed to the quenched CLT; the proof is motivated by ideas in \cite[\S6-7]{\PZ} and \cite[Lem.~4.1]{\BolthausenSznitman}. For given $n\ge1$, let $\bs$ denote the unique increasing bijection
$$\bs:\Z_{\ge0} \to \bigcup_{i\ge1}[\tau_i,\eta_i),$$
and let $X^\exc_t\equiv X_{\bs(t)}$, the excursion process of $X$ with parameter $n$ (recalling Rmk.~\ref{r:param.n}). For $\bs(t)\in[\tau_i,\eta_i)$ write $X^\cent_t\equiv \abs{X^\exc_t}-\abs{X_{\tau_i}}$.

\bpf[Proof of Thm.~\ref{t:mclt}]
We show the quenched CLT for $X$ through a quenched CLT for $X^\cent$ along geometrically increasing subsequences $b_k\equiv \flr{b^k}$ ($k\ge0$) with $b>1$.

\medskip\noindent
{\it Step 1: annealed CLT for $X^\cent$.} \\
The time killed during the first $n$ steps of $X$ is $n-\bs^{-1}(n) \le \De_n$, so Cor.~\ref{c:azuma}~(\ref{c:azuma.int}) gives $n^{-1} \sup_{0\le t\le T} \abs{\bs(\flr{nt})-\flr{nt}}\to0$ in $\mgw$-probability. It follows from Propn.~\ref{p:mgw.aclt} and the continuity of Brownian motion that the processes $X^\cent_{\flr{nt}}/(\si\sqrt{n})$ also satisfy the annealed $\mgw$-CLT.

\medskip\noindent
{\it Step 2: quenched CLT for $X^\cent$ along geometrically increasing subsequences.} \\
Recalling Rmk.~\ref{r:polyg}, let $\B^n(X)\equiv(\B^n_t(X))_{t\ge0}$ denote the polygonal interpolation of $j/n\mapsto X^\cent_j/(\si\sqrt{n})$, and regard $\B^n(X)$ as an element of $C[0,T]$ with the norm
$$d_T(u,u')
\equiv \Big(
	\sup_{0\le t\le T} \abs{u_t-u'_t}
	\Big) \wedge 1.$$
We will show that for all Lipschitz functions $F:C[0,T]\to[-1,1]$ with Lipschitz constant $\le1$,
\beq\label{e:bs.var}
\sum_{k\ge0}
	\Var_\mgw[\E_\tree[F[\B^{\flr{b^k}}(X)] ]]
< \infty.
\eeq
The Borel--Cantelli lemma then implies (cf.\ \cite[Lem.~4.1]{\BolthausenSznitman}) that for $\mgw$-a.e.\ $\tree$, the processes $X^\cent_{\flr{nt}}/(\si\sqrt{n})$ converge in law to the absolute value of a standard Brownian motion along the subsequence $b_k$.

To see \eqref{e:bs.var}, let $\tree\sim\mgw$, let $(X^i,\bs^i)$ be two independent realizations of $(X,\bs)$ conditioned on $\tree$, and write $\B^{n,i}\equiv\B^n(X^i)$. Then
$$\Var_\mgw[\E_\tree[F[\B^{n}(X)] ]]
= \E_\mgw[ F(\B^{n,1}) F(\B^{n,2}) ]
	-\E_\mgw[ F(\B^{n,1})]^2.$$
Let $\cE^i_n$ denote the subtree explored by $X^i$ up to time $n$. Conditioning on the first $\ell(n)/2$ levels of $\tree$, let $(\acute\cE^i_n,\acute X^i,\acute\bs^i)$ ($i=1,2$) be two independent realizations of $(\cE^i_n,X^i,\bs^i)$: then the processes $(\acute X^i)^\cent$ are \emph{exactly} independent with law not depending on the first $\ell(n)/2$ levels of $\tree$. Moreover, if $\cA_n$ denotes the event that the paths of $\acute X^1$ and $\acute X^2$ up to time $\max_i (\acute\bs^i)^{-1}(n)$ have no common vertices at distance more than $\ell(n)/2$ from the root, then we can couple $(\cE^i_n,X^i|_{[0,n]})_{i=1,2}$ with $(\acute\cE^i_n,\acute X^i|_{[0,n]})_{i=1,2}$ such that the processes agree on the event $\cA_n$. Therefore
\begin{align*}
&\Var_\mgw[\E_\tree[F[\B^{n}(X)] ]]
\le \E_\mgw[ F(\acute\B^{n,1}) F(\acute\B^{n,2}) ]
	+ \mgw(\cA_n)
	-\E_\mgw[ F(\B^{n,1})]^2\\
&=\E_\mgw[ F(\acute\B^{n,1})]^2+\mgw(\cA_n)
	-\E_\mgw[ F(\B^{n,1})]^2
=\mgw(\cA_n)
\end{align*}
We claim $\mgw(\cA_n)\le n^{-c}$: since Cor.~\ref{c:timeroot} and Cor.~\ref{c:azuma}~(\ref{c:azuma.int}) imply $\mgw(2n-\bs^{-1}(2n) \ge n) \le n^{-c}$, it suffices to bound the probability that the paths of $X^1$ and $X^2$ up to time $2n$ intersect at distance $>\ell(n)/2$ from the root. But the chance that $X^2$ hits a given vertex $v$ with $\abs{v}>\ell(n)/2$ by time $2n$ is $\le Cn\rho^{-\ell(n)/2}$, and summing over the vertices visited by $X^1$ proves the claim. The variance condition \eqref{e:bs.var} now follows by summing over $(b_k)_{k\ge0}$.

\medskip\noindent
{\it Step 3: quenched CLT for $X$ along geometrically increasing subsequences.} Extend $\bs^{-1}$ to a nondecreasing map $\Z_{\ge0}\to\Z_{\ge0}$ by setting $\bs^{-1}(t)=\bs^{-1}(\tau_i)$ for $t\in[\eta_{i-1},\tau_i)$: then
$$\abs{\abs{X_t}-X^\cent_{\bs^{-1}(t)}}
=\begin{cases}
\abs{X_{\tau_i}}, & t\in[\tau_i,\eta_i), \\
\abs{X_t}, & t\in[\eta_{i-1},\tau_i).
\end{cases}$$
It follows from Lem.~\ref{l:azuma} and Cor.~\ref{c:azuma}~(\ref{c:azuma.int}) that for any $b>1$,
$$\left.\begin{array}{r}
b_k^{-1}
	\sup_{0\le t\le b_k T}
	\{
	\bs^{-1}(t)-\bs^{-1}(t)
	\}\\
b_k^{-1/2}
	\sup_{0\le t\le b_kT}
	\abs{\abs{X_t}-\abs{X^\cent_{\bs^{-1}(t)}}}
\end{array}\right\}\stackrel{k\to\infty}{\longrightarrow}0,\quad
	\mgw\text{-a.s.}$$
It follows that the processes $\abs{X_{nt}}/(\si\sqrt{n})$ satisfy the quenched CLT along the subsequence $(b_k)_{k\ge0}$ for any $b>1$.

\medskip\noindent
{\it Step 4: quenched CLT for $X$ along full sequence.} For the processes $(\abs{X_{\flr{nt}}}/(\si\sqrt{n}))_{t\ge0}$, $\mgw$-a.s.\ tightness and convergence of finite-dimensional distributions both follow from the scaling relation
$$\B^n_t(X)
=\sqrt{\f{b_k}{n}}
\B^{b_k}\Big(t \f{n}{b_k}\Big),$$
(cf.\ proof of \cite[Lem.~4.1]{\BolthausenSznitman}).
\epf

\bpf[Proof of Cor.~\ref{c:mgw.clt.cts}]
Given $(\tree,X)\sim\mgw\otimes\rw_\rho$ we can obtain $(\tree,X^\cts)\sim\mgw\otimes\rw_\rho^\cts$ by taking $(E_i)_{i\ge1}$ i.i.d.\ exponential random variables with unit mean independent of $X$, and setting
$$X^\cts_t = X_{\bthet(t)},\quad
\bthet(t) = \max\setb{i : \sum_{j=1}^i \f{E_j}{\rho+d_{X_{j-1}}} \le t};$$
similarly we can obtain $((\tree^\coup,\ray),Y^\cts)\sim\imgw_0\otimes\rw_\rho^\cts$ from $((\tree^\coup,\ray),Y)\sim\imgw_0\otimes\rw_\rho$. Thus a shifted coupling of $(\tree,X)$ with $((\tree^\coup,\ray),Y)$ (as constructed in \S\ref{s:coup}) naturally gives rise to a shifted coupling of $(\tree,X^\cts)\sim\mgw\otimes\rw_\rho^\cts$ with $((\tree^\coup,\ray),Y^\cts)\sim\imgw_0\otimes\rw_\rho^\cts$ by using sequences $(E_i)_{i\ge1}$ for $X^\cts$ and $(E^\coup_i)_{i\ge1}$ for $Y^\cts$ which are marginally i.i.d.\ exponential but such that the jump times match during the coupled excursions.

By Thm.~\ref{t:imgw}~(\ref{t:imgw.erg}) and the exponential decay of the $E^\coup_i$, it holds $\imgwr$-a.s.\ that
$$\f{1}{n}\sum_{i=1}^n \f{E^\coup_i}{\rho+d_{Y_i}} \to
\E_\imgwr\Big[\f{1}{\rho+d_\rt}\Big]
= \f{1}{2\rho}.$$
From this it is easy to see that $n^{-1} \sup_{0\le t\le T}[\bthet(nt)-2\rho nt]\to0$ $\imgwr$-a.s., so on $\imgwr$-a.e.\ $(\tree,\ray)$ the processes $(h(Y^\cts_{\flr{nt}})/(\si\sqrt{2\rho n}))_{t\ge0}$ converge in law to standard Brownian motion. The quenched $\mgw$-CLT for $X^\cts$ follows from the proof of Thm.~\ref{t:mclt}.
\epf

\section{Transience-recurrence boundary for \texorpdfstring{$\rwre_\lm$}{RWRE(lambda)}}
\label{s:rwre.trans.recur}

We now prove Thm.~\ref{t:tr}. Our proof is a straightforward adaptation of that of \cite[Thm.~1]{\LPRWRE} or \cite[Propn.~1.1]{\Faraud} once we supply the needed large deviations estimate (Lem.~\ref{l:rwre.large.devs}) on the conductances at the $n$-th level of the tree, extending the estimates of \cite[p.~129]{\LPRWRE} and \cite[p.~7]{\Faraud} to our setting of Markovian dependency.

Let $\domain \equiv \set{\gam : \max_{a,b}\bar A^{(\gam)}(a,b)<\infty}$, where $\bar A^{(\gam)}$ is as defined in \eqref{e:rwre.mat}. Recall that $\bar\rho(\gam)$ denotes the Perron--Frobenius eigenvalue of $\bar A^{(\gam)}$ with $\bar\rho(\gam)\equiv\infty$ for $\gam\notin\domain$. The following lemma collects some basic properties of $\bar\rho$.

\blem\label{l:bar.rho}
Under the hypotheses of Thm.~\ref{t:tr}, $\bar\rho$ is lower semi-continuous and log-convex on $\R$, and differentiable on $\domain$.

\bpf
Lower semi-continuity of $\bar\rho$ in $\domain$ follows from Fatou's lemma, and lower semi-continuity outside the closure of $\domain$ is trivial, so it remains to consider the boundary of $\domain$: we must show that if $\gam\to\gam_\infty$ with $\max_{a,b} \bar A^{(\gam)}(a,b)\to\infty$ then $\bar\rho(\gam)\to\infty$. Recall the min-max characterization (see e.g.\ \cite[Cor.~8.3.3]{\HJ})
$$\bar\rho(\gam)
= \max_{\vec x\ge0, \vec x\ne0}
\min_{a : x_a\ne0} \f{(\bar A^{(\gam)} \vec x)_a}{x_a}.$$
Since $\rho[(A^{(\gam)})^k]=\rho(A^{(\gam)})^k$ ($k\in\N$), and $\bar A^{(0)}$ is positive regular which implies $\bar A^{(\gam)}$ is also for all $\gam\in\domain$, we may assume without loss that $\min_{a,b} \bar A^{(\gam)}(a,b)\ge\ep$ for all $\gam$ in a neighborhood of $\gam_\infty$. Applying the min-max characterization to the vectors $\vec x=\I_a$ gives $\bar\rho(\gam) \ge \max_a \bar A^{(\gam)}(a,a)$. Applying it to the vectors
$$\vec x
=\Big( \Ind{c=a} + \Ind{c=b}
	\f{\sqrt{\ep}}{ \sqrt{A^{(\gam)}(a,b)}} \Big)_{c\in\cQ},
	\quad a\ne b$$
gives
\begin{align*}
\bar\rho(\gam)
&\ge
\Big(
{\textstyle A^{(\gam)}(a,a)+\sqrt{\ep A^{(\gam)}(a,b)}}
\Big)
\wedge
\Big(
A^{(\gam)}(b,a) \f{\sqrt{A^{(\gam)}(a,b)}}{\sqrt{\ep}}
+ A^{(\gam)}(b,b)
\Big)\\
&\ge \sqrt{\ep A^{(\gam)(a,b)}}.
\end{align*}
Combining gives $\bar\rho(\gam)^2 \ge \ep\max_{a,b\in\cQ} A^{(\gam)}(a,b)$ which proves lower semi-continuity.

The entries of $\bar A^{(\gam)}$ are log-convex in $\gam$ by H\"older's inequality, so $\bar\rho$ is log-convex by monotonicity and log-convexity of the Perron--Frobenius eigenvalue in the entries of the matrix
(see e.g.\ \cite[Cor.~8.1.19]{\HJ} and \cite[Exercise 4.34]{\BVCvx}). For differentiability of $\bar\rho$ in $\domain$ see \cite[p.~75]{\DZ}.
\epf
\elem

For $\gam\in\domain$ let $\rvre^{(\gam)}$ and $\lvre^{(\gam)}$ denote the associated left and right Perron--Frobenius eigenvectors; we use the shorthand
$$\rv\equiv \rvre^{(0)},\quad\lv\equiv\lvre^{(0)}.$$
For $\tree\sim\mgwre$ and $v\in\tree$, let
$$C_v \equiv \prod_{\rt<u\le v} \al_u,$$
the conductance of the edge leading to $v$. The natural generalization of the martingale introduced in \S\ref{ss:imgwr.wk.lim} is
\beq\label{e:rwre.mg}
\Znorm_n^{(\gam)} = \f{1}{\bar\rho(\gam)^n} \sum_{v\in D_n} \rvecre{\type_v}^{(\gam)} C_v^\gam;
\eeq
this is a multi-type Mandelbrot's martingale and has been studied in various contexts, for example as the Laplace transform of the branching random walk with increments $\log\al_v$ \cite{\BNT,\KSBRW}. Using this martingale we can make a change of measure and control the conductances at the $n$-th level by controlling the conductance of the edge leading to a \emph{random} vertex: recalling \eqref{e:Qn.a}, for each $a\in\cQ$ define the size-biased measure $\Qnre^a_n$ on $\trees$ by
$$\f{d\Qnre^a_n}{d \mgwretype{a}}
=\f{\Znorm^{(0)}_n}{\rvec{a}}.$$
We then let $\Qnre^a_{n\star}$ denote the measure on pairs $(\tree,v_n)$ obtained by letting $\tree\sim\Qnre^a_n$ and choosing $v_n\in D_n$ according to weights $\rvec{\type_v}$.

\blem\label{l:rwre.large.devs}
Under the hypotheses of Thm.~\ref{t:tr}, for each $a\in\cQ$, under $\Qnre^a_{n\star}$ the random variables $n^{-1} \log C_{v_n}$ satisfy a large deviation principle with good rate function
$\Lm^*(x) \equiv \sup_\gam (\gam x - \Lm(\gam))$,
where $\Lm(\gam)\equiv\log\bar\rho(\gam)-\log\bar\rho(0)$. In particular, for any $0<z<y$,
\beq \label{e:rwre.large.devs}
\liminf_{n\to\infty} \f{1}{n} \log \Qnre^a_{n\star} (C_{v_n} > z^n)
\ge - \sup_{\gam\ge0} (\gam \log y - \Lm(\gam)).
\eeq

\bpf
Fixing $a\in\cQ$, let $\Lm_n\equiv\Lm^a_n$ denote the cumulant generating function of $n^{-1} \log C_{v_n}$ with respect to $\Qnre^a_{n\star}$, that is,
$$\Lm_n(\gam) = \log \E_{\Qnre^a_{n\star}}[C_{v_n}^{\gam/n}].$$
Then
\begin{align*}
e^{\Lm_n(n\gam)}
&= \E_{\Qnre^a_{n\star}}[C_{v_n}^\gam]
= \E_{\Qnre^a_{n}}\bigg[
	\f{ \sum_{v\in D_n} \rvec{\type_v} C_v^\gam }
	{\sum_{v\in D_n} \rvec{\type_v}} \bigg]
= \f{1}{\rvec{a} \bar\rho(0)^n}
	\E_{\mgwretype{a}}
	\bigg[ \sum_{v\in D_n} \rvec{\type_v}  C_v^\gam \bigg] \\
&\asymp \f{\bar\rho(\gam)^n}{\bar\rho(0)^n}
	\E_{\mgwretype{a}}[\Znorm^{(\gam)}_n]
\asymp \f{\bar\rho(\gam)^n}{\bar\rho(0)^n},
\end{align*}
where $\asymp$ indicates equivalence up to constant factors depending only on $\rv$ and $\rvre^{(\gam)}$. Thus
$$\lim_{n\to\infty} \f{1}{n} \Lm_n(n\gam)
= \log\bar\rho(\gam) - \log\bar\rho(0) = \Lm(\gam).$$
By Lem.~\ref{l:bar.rho} this is an essentially smooth convex function in the sense of \cite[Defn.~2.3.5]{\DZ}, so the large deviation principle follows from the G\"artner-Ellis theorem (see \cite[Thm.~2.3.6]{\DZ}). In particular, for any $0<z<y$, \cite[(2.3.8)]{\DZ} implies
$$\liminf_{n\to\infty} \f{1}{n} \log \Qnre^a_{n\star}(C_{v_n} > z^n)
\ge -\inf_{x > \log z} \Lm^*(x) \ge -\Lm^*(\log y)$$
(making use of \cite[Lem.~2.3.9]{\DZ}). The result \eqref{e:rwre.large.devs} follows immediately if $\log y=\Lm'(\gam)$ for some $\gam\ge0$, or if $\log y\ge\sup_{\gam\ge0}\Lm'(\gam)$ in which case $\sup_\gam(\gam\log y-\Lm(y))=\lim_{\gam\to\infty}(\gam\log y-\Lm(y))$. Next, the assumption that $\Lm<\infty$ in a neighborhood of $0$ implies, via the relation $\Lm(\gam)=\sup_x(x\gam-\Lm^*(x))$, that $\lim_{\abs{x}\to\infty}\Lm^*(x)=\infty$, therefore $\Lm^*$ attains its global infimum at $x_0=\Lm'(0)$ with $\Lm^*(x_0)=-\Lm(0)$. Therefore \eqref{e:rwre.large.devs} again holds in the remaining case $\log y\le\inf_{\gam\ge0}\Lm'(\gam)=x_0$.
\epf
\elem

Thm.~\ref{t:tr} now follows by adapting the proof of \cite[Thm.~1]{\LPRWRE}:

\bpf[Proof of Thm.~\ref{t:tr}]
Since the bias $\lm$ can always be absorbed into the environment variables $\al_v$ ($v\in\tree$), we may take $\lm=1$ from now on, and write
$p \equiv p_1 = \min_{0\le \gam\le 1} \bar\rho(\gam)$.

\medskip\noindent
(\ref{t:tr.r}) Suppose $p<1$. We will use the fact that the random walk is positive recurrent if and only if the conductances have finite sum \cite[Propn.~9-131]{\KSK}. If $\bar\rho(\gam)<1$ for some $\gam\in[0,1]$ then
$$\E\Big[ \sum_{v\in\tree} C_v^\gam \Big]
\asymp \sum_{n\ge0} \bar\rho(\gam)^n < \infty,$$
so $\sum_{v\in\tree} C_v^\gam<\infty$ a.s. In particular $C_v\le C_v^\gam<1$ for all but finitely many $v\in\tree$ so $\sum_{v\in\tree} C_v<\infty$ a.s.

\medskip\noindent
(\ref{t:tr.t}) Suppose $p>1$. We will show that on the event of non-extinction there exists $w<1$ such that
\beq \label{e:rwre.exp.decrease}
\liminf_{n\to\infty} w^n \sum_{v\in D_n}C_v > 0;
\eeq
transience then follows from \cite[Cor.~4.2]{\Lyons}. By the proof on \cite[p.~129]{\LPRWRE},
$$p = \max_{0<y\le 1}
	\Big\{
	y \, \inf_{\gam\ge0} y^{-\gam} \bar\rho(\gam)
	\Big\},$$
and we fix $y\in(0,1]$ achieving this maximum. Then Lem.~\ref{l:rwre.large.devs} implies that
$$\liminf_{n\to\infty} \f{1}{n} \log \Qnre^a_{n\star}(C_{v_n}>z^n)
\ge - \log [y\bar\rho(0)/p]\quad \forall z<y, \forall a\in\cQ.$$
Therefore we can choose $z<y$, $\ell\in\N$, and $\ep,w\in(0,1)$ such that
$$\min_{a\in\cQ} \Qnre^a_{\ell\star}(
	\set{C_{v_\ell}>z^\ell}
	\cap \set{\al_w\ge\ep \ \forall \rt< w \le v_\ell})
	\ge q > (w z \bar\rho(0))^{-\ell}.$$
Now consider the following percolation process (same as on \cite[p.~130]{\LPRWRE}): let $\tree\sim\mgw(\cdot\giv\nx^c)$, and let $\tree[\ell]$ be the tree with vertices $\set{v\in\tree : \abs{v}\equiv0 \text{ mod } \ell}$, with an edge $v\to w$ if and only if $\abs{w}=\abs{v}+\ell$ in $\tree$. Form a random subgraph $\tree[\ell]^\perc\subseteq\tree[\ell]$ by keeping the edge $v\to w$ if and only if
$$\prod_{v<u\le w} \al_u > z^\ell \quad\text{ and } \quad
\min_{v<u\le w} \al_u \ge\ep,$$
in which case we write $v\rightsquigarrow w$. The subtree of $\tree[\ell]^\perc$ descended from any vertex $v$ has the law of a multi-type Galton--Watson tree with mean offspring numbers
$$A^\perc(a,b)
= \E_{\mgwretype{a}}
	\bigg[ \sum_{v\in D_\ell} \Ind{\type_v=b} \Ind{\rt\rightsquigarrow v} \bigg]
= \f{\rvec{a} \bar\rho(0)^\ell}{\rvec{b}}
	\E_{\Qnre^a_{\ell\star}} [\Ind{\type_{v_\ell}=b} \Ind{\rt\rightsquigarrow v_\ell} ].$$
We calculate
$$\sum_b A^\perc(a,b) \rvec{b}
\ge \rvec{a} \bar\rho(0)^\ell q > \rvec{a} (wz)^{-\ell},$$
so $A^\perc$ has Perron--Frobenius eigenvalue larger than $(wz)^{-\ell}$, and consequently $\tree[\ell]^\perc$ a.s.\ has a connected component which is an infinite tree $\tree[\ell]^\star$ of branching number larger than $(wz)^{-\ell}$, rooted at some $o^\star\in\tree[\ell]^\perc$. It follows that the left-hand side of \eqref{e:rwre.exp.decrease} is
$$\ge C_{\rt^\star} w^{\abs{\rt^\star}\ell}
	\ep^{\ell-1}
	\liminf_{n\to\infty}
	\sum_{v\in D_n(\tree[\ell]^\star)}
	(wz)^{n\ell}>0$$
which concludes the proof.
\epf

\section{Appendix: general properties of \texorpdfstring{$\mgw$}{MGW} trees}
\label{s:appx}

In this section we prove some basic facts about $\mgw$ trees which were used in the proof of the main theorem. In \S\ref{ss:appx.mmt} we prove Propn.~\ref{p:moments} which states that \hypmmt{\mmt} implies $\E_\mgw[W_\rt^\mmt]<\infty$. In \S\ref{ss:appx.cond} we prove a conductance lower bound (Propn.~\ref{p:polar}) which gives Lem.~\ref{l:polar}. We begin in \S\ref{ss:appx.gen} by collecting some preliminary observations.

\subsection{Generating function and subtree of infinite descent}\label{ss:appx.gen} Let
$$F(\vec s)\equiv (F^a(\vec s))_{a\in\cQ}
\equiv \Big(
	\E_{\q^a}\Big[ \prod_{b\in\cQ} s_b^{x_b} \Big]
	\Big)_{a\in\cQ},\quad
	\vec s\in[0,1]^\cQ;$$
we refer to $F$ as the generating function of the $\mgw$ tree. If $F^{(n)}$ denotes the $n$-fold composition of $F$, then for all $a\in\cQ$
$$\E_{\mgwtype{a}}
	\bigg[ \prod_{b\in\cQ} s_b^{Z_n(b)}\bigg]
= (F^{(n)}(\vec s))_a,\quad
\mgwtype{a}(\abs{\vec Z_n}=0) = (F^{(n)}(\vec 0))_a.$$
Next let
$$\Phi(s)\equiv(\Phi^a(s))_{a\in\cQ}
\equiv\Big(
	\E_{\mgwtype{a}}[e^{-s W_\rt}]
	\Big)_{a\in\cQ},\quad
s\ge0,$$
and let $\phi(s)\equiv\ip{\lv}{\Phi(s)}=\E_\mgw[e^{-sW_\rt}]$. Since
$$\Phi(s)=\lim_{n\to\infty} F^{(n)}
	((\exp\{-s\rvec{b} /\rho^n\})_{b\in\cQ}),$$
we have the functional relation $\Phi(s)=F[\Phi(s/\rho)]$.

For many purposes the case of $\mgw(\nx)\in(0,1)$ can be reduced to the simpler case of an a.s.\ infinite tree without leaves by the following transformation which is discussed in \cite[\S I.12]{\AthreyaNey} for the single-type case. For $\tree\sim\mgw$, consider the subtree $\tree^\infty$ consisting of those vertices $v$ of \emph{infinite descent}, i.e.\ with $\abs{\subtree{\tree}{v}}=\infty$. Conditioned on $\nx^c$, $\tree^\infty$ is an a.s.\ infinite tree without leaves, following a transformation of the original $\mgw$ given by generating function $\acute F(\vec s) = (\acute F^a(\vec s))_{a\in\cQ}$, where, with $\mathfrak{x}_a\equiv \mgwtype{a}(\nx)$,
\begin{align*}
\acute F^a(\vec s)
&\equiv \f{1}{1-\mathfrak{x}_a} \sum_{\vec x} \q^a(\vec x) \prod_{b\in\cQ}
	\sum_{y_b\le x_b} {x_b \choose y_b} \mathfrak{x}_b^{x_b-y_b} (1-\mathfrak{x}_b)^{y_b} s_b^{y_b}\\
&= \f{F_a((\mathfrak{x}_b + (1-\mathfrak{x}_b)s_b)_{b\in\cQ})}{1-\mathfrak{x}_a}.
\end{align*}
The transformed law has mean matrix
$$\acute A = D^{-1} AD,\quad D = \diag( (1-\mathfrak{x}_a)_{a\in\cQ} ),$$
so in particular it has the same Perron--Frobenius eigenvalue as $A$. Finally, it is clear that if the original law satisfies \hypmmt{\mmt} then so does the transformed law.

\subsection{Positive moments of the normalized population size}
\label{ss:appx.mmt}

In this section we show that moment conditions on the $\mgw$ offspring distribution translate directly to moment conditions on the normalized population size of the entire tree. We begin by recalling an easy fact concerning Laplace transforms (see e.g.\ \cite[\S XIII]{\FellerII}).

\blem\label{l:laplace}
Let $\vph(s)\equiv\E[e^{-sW}]$ be the Laplace transform of a non-negative random variable $W$. For any integer $n\ge0$, $\E[W^n]<\infty$ if and only if there exist finite coefficients $m_0,\ldots,m_n$ such that
\beq\label{e:lapl.rem}
\sum_{r=0}^n \f{m_r}{r!} s^r
=\vph(s)+ o(s^n),\quad s\decto0.
\eeq
In this case $m_r=\E[(-W)^r]=\lim_{s\decto0}\vph^{(r)}(s)$, and the left-hand side of \eqref{e:lapl.rem} is the $n$-th order (one-sided) Taylor expansion $\poly_{n,0}\vph$ of $\vph$ at $0$.

\bpf
($\Ra$) The function $\vph$ is infinitely differentiable on $(0,\infty)$ with $n$-th derivative given by $(-1)^n\vph^{(n)}(s)=\E[e^{-sW} W^n]$, and by the monotone convergence theorem
$$\vph^{(n)}(0)\equiv\lim_{s\decto0} \vph^{(n)}(s)
=\E[W^n]\in[0,\infty].$$
Writing $\expm(x)\equiv e^{-x}$, by Taylor's theorem
$$\rem_{n,0} \expm(x)
\equiv (-1)^{n+1}[e^{-x}-\poly_{n,0} \expm(x)]
=  \f{x^n}{n!} (1-e^{-\ze}),\quad
0\le\ze\le x.$$
If $\E[W^n]<\infty$ then
$$\E[\rem_{n,0}\expm(sW)]
= \f{s^n}{n!} \E[ W^n(1-e^{-\ze W}) ],
\quad 0\le\ze\le s.$$
The right-hand side is $o(s^n)$ so \eqref{e:lapl.rem} holds with
\beq\label{e:lapl.rem.expm}
\vph(s)-\sum_{r=0}^n \f{m_r}{r!} s^r
=\vph(s)- \poly_{n,0}\vph(s)
= (-1)^{n+1}\E[\rem_{n,0}\expm(sW)].
\eeq

\medskip\noindent
($\La$) Assuming \eqref{e:lapl.rem}, suppose inductively that $m_r=\vph^{(r)}(0)$ for $0\le r\le k$ with $k<n$. For $0<s_0\le s$ we have
$$\vph(s)
= \poly_{k,s_0}\vph(s)
	+ \f{\vph^{(k+1)}(\ze)}{(k+1)!}(s-s_0)^{k+1},\quad
	s_0 \le \ze\le s,$$
so by \eqref{e:lapl.rem} and the inductive hypothesis
$$
o(s^{k+1})
=\vph(s)-\sum_{r=0}^{k+1} \f{m_r}{r!} s^r
=\sum_{r=0}^k \f{o(s_0)}{r!}
	+\f{\vph^{(k+1)}(\ze)(s-s_0)^{k+1}
	-m_{k+1}s^{k+1}}{(k+1)!}.$$
Taking $s_0 \ll s^{k+1}$ we find a contradiction unless $\lim_{\ze\decto0}\vph^{(k+1)}(\ze)=m_{k+1}$.
\epf
\elem

\blem\label{l:moments.int}
If \hypirred, \hypks\ and \hypmmt{n} hold with $n\in\Z_{\ge2}$, then $\E_\mgw[W_\rt^n]<\infty$.

\bpf
Following the proof of \cite[Thm.~0]{\BinghamDoney}, we will show $\E_\mgw[W_\rt^n]<\infty$ using the characterization Lem.~\ref{l:laplace} of the derivatives at zero of the Laplace transform $\phi(s)=\E_\mgw[e^{-sW_\rt}]$. Write $\bS_\cQ\equiv\set{\vec v\in[0,\infty)^\cQ:\sum_{a\in\cQ} v_a=1}$, and define
$$f(t;\vec v)
\equiv\ip{\lv}{F(e^{-t\vec v})}
=\E_\mgw[e^{-t\ip{\vec v}{\vec Z_1}}],\quad
t\ge0,\,\vec v\in\bS_{\cQ}.$$
By \eqref{e:lapl.rem.expm},
$$f(t;\vec v)=\poly_{n,0}f(t;\vec v)+(-1)^{n+1}\rem_{n,0}f(t;\vec v),\quad
\lim_{t\decto0} \Big( t^{-n}
	\sup_{\vec v\in\bS_{\cQ}}
	\rem_{n,0}f(t;\vec v)
	\Big)=0,$$
where $\poly_{n,0}f(t;\vec v)$ is a polynomial of degree at most $n$ in the entries of $t\vec v$ satisfying
\beq\label{e:poly.n}
\poly_{n,0}f(t;\vec v)
=1-t\ip{\vec v}{\E_\mgw[\vec Z_1]}+O(t^2)
=1-\rho t\ip{\vec v}{\lv}+O(t^2).
\eeq
If we let $t\equiv t(s)\ge0$ and $\vec v\equiv \vec v(s)\in\bS_{\cQ}$ be defined by $\Phi(s/\rho)=e^{-t\vec v}$, then
\beq\label{e:phi.s.t}
\phi(s)
=\ip{\lv}{\Phi(s)}
=\ip{\lv}{F[\Phi(s/\rho)]} = f(t;\vec v) = \poly_{n,0}f(t;\vec v)+o(t^n)
\eeq
(using that $\rem_{n,0}f(t;\vec v)=o(t^n)$ uniformly over $\vec v\in\bS_{\cQ}$).

We next expand $t\vec v$ in powers of $s$. Note that $\phi(s)-m_0-m_1 s=o(s)$ where $m_0=1$, $m_1=-\E_\mgw[W_\rt]$, so suppose inductively that $\E_\mgw[W_\rt^k]<\infty$ for some $1<k<n$. Lem.~\ref{l:laplace} implies the existence of polynomials $q^a_k(s)=\rvec{a}+O(s)$ such that
\beq\label{e:t.in.s}
tv_a
= -\log\Phi^a(s/\rho)
= s\,q^a_k(s) + o(s^k),\quad s\decto0.
\eeq
By recalling \eqref{e:phi.s.t} and comparing \eqref{e:poly.n} against the Taylor expansion of $\rho\ip{\lv}{e^{-t\vec v}}$, we find
$$[\phi(s)-1]-\rho[\phi(s/\rho)-1]
= Q_n(t\vec v) + o(t^n),\quad t\decto0.$$
for $Q_n:\R^\cQ\to\R$ polynomial with $Q_n(t\vec v)=O(t^2)$. But squaring \eqref{e:t.in.s} gives $(tv_a)^2 = s^2 q^a_k(s)^2 + o(s^{k+1})$, and substituting into the above and dividing through by $s$ gives that
$$\psi(s)-\psi(s/\rho)
\equiv\f{\phi(s)-1}{s}-\f{\phi(s/\rho)-1}{s/\rho}
= s P_k(s) + o(s^k)$$
for $P_k$ polynomial in $s$. Since $\rho>1$ and $\lim_{s\decto0}\psi(s)=\phi'(0)$,
$$\psi(s) = \psi'(0) + \sum_{j\ge0} (s/\rho^j)P_n(s/\rho) + o(s^k) = \wt P_n(s) + o(s^k)$$
for $\wt P_n$ another polynomial in $s$. This verifies the inductive hypothesis by the definition of $\psi$ together with another application of Lem.~\ref{l:laplace}.
\epf
\elem

If \hypmmt{n} holds with $n\in\Z_{\ge2}$ then we write $\phi(s)=\poly_{n,0}\phi(s)+(-1)^{n+1}\rem_{n,0}\phi(s)$ with $\poly_{n,0}\phi$ polynomial of degree at most $n$ and $\rem_{n,0}\phi(s)=o(s^n)$. An easy consequence of the proof of the Lem.~\ref{l:moments.int} is the following

\bcor\label{c:moments.int}
If \hypmmt{n} holds with $n\in\Z_{\ge2}$ then
$$\rem_{n,0}\phi(s)-\rho\rem_{n,0}\phi(s/\rho)
=\rem_{n,0}f(t,\vec v)+O(t^{n+1})$$
for $\Phi(s/\rho)\equiv e^{-t\vec v}$.

\bpf
Summing \eqref{e:t.in.s} over $a\in\cQ$ gives the existence of $b_2,\ldots,b_n$ finite such that
$$t = s + \sum_{r=2}^n b_r s^r + o(s^n).$$
It follows easily that $s$ has a similar expansion in terms of $t$: indeed $s=t+o(t)$, so suppose inductively that for some $1\le k<n$ there exist $c_2,\ldots,c_k$ finite such that $s=t + \sum_{r=2}^k c_r t^r+o(t^k)$. Then
$$s=t-\sum_{r=2}^n b_r
	\Big(t + \sum_{r=2}^k c_r t^r+o(t^k) \Big)^2
	+ o(s^n),$$
which is a polynomial in $t$ plus $o(t^{k+1})$. This verifies the inductive hypothesis so we conclude that $s=t+\sum_{r=2}^n c_r t^r + o(t^n)$ as claimed. From the proof of Lem.~\ref{l:moments.int} we have
\begin{align*}
o(s^n)&=(-1)^{n+1}
[\rem_{n,0}\phi(s)-\rho\rem_{n,0}\phi(s/\rho)]\\
&=f(t;\vec v)
-\rho\ip{\lv}{e^{-t\vec v}}
-[\poly_{n,0}\phi(s)
-\rho \poly_{n,0}\phi(s/\rho)]\\
&=f(t;\vec v) - Q_n^\backprime(t\vec v)- s^2 q(s) +O(t^{n+1})
\end{align*}
for $Q_n^\backprime$ and $q$ polynomial. But by the above $s^2$ can be expressed as a polynomial in $t$ up to $o(t^{n+1})$ error, so in fact
$$o(t^n)=(-1)^{n+1}
[\rem_{n,0}\phi(s)-\rho\rem_{n,0}\phi(s/\rho)]
=f(t;\vec v)-Q_n^{\backprime\backprime}(t\vec v)+O(t^{n+1})$$
for $Q_n^{\backprime\backprime}:\R^\cQ\to\R$ polynomial in $t\vec v$ of degree at most $n$ in $t$, whence necessarily $Q_n^{\backprime\backprime}(t\vec v)=\poly_{n,0}f(t;\vec v)$ as claimed.
\epf
\ecor

If $\Phi(s/\rho)=e^{-t\vec v}$ with $t\equiv t(s)\ge0$, $\vec v\equiv \vec v(s)\in\bS_{\cQ}$ as above, then
$$t=-\sum_{a\in\cQ}\log\E_{\mgwtype{a}}[e^{-sW_\rt}].$$
In particular $t'(s)$ is finite and positive for all $s>0$ with $\lim_{s\decto0}t'(s)=1$, so $s\mapsto t(s)$ is an increasing bijection from $[0,\infty)$ to $[0,t_{\max})$ where $t_{\max}=-\sum_{a\in\cQ}\log \mgwtype{a}(\nx)$. For $t<t_{\max}$ we therefore write $\vec v_t\equiv \vec v(s)$ with $s$ defined by $t=t(s)$.

\bpf[Proof of Propn.~\ref{p:moments}]
Since the subtree of infinite descent described in \S\ref{ss:appx.gen} has the same normalized population size as the original tree, we may reduce to the case $\mgw(\nx)=0$ so that $t_{\max}=\infty$.

By Lem.~\ref{l:moments.int} we may take $\mmt=n+\be$ for $\be\in(0,1)$, and by \cite[Thm.~B]{\BinghamDoney} the result follows upon showing
$$\int_0^1 \rem_{n,0}\phi(s) \, s^{-(1+\mmt)} \d s<\infty.$$
By the proof of \cite[Propn.~5]{\BinghamDoney}, this in turn follows upon showing
\beq\label{e:mgw.laplace.integral}
\int_0^1 \rem_{n,0}f(t;\vec v_t) \, t^{-(1+p)}\d t<\infty.
\eeq
(replacing \cite[(3.13)]{\BinghamDoney} with \eqref{e:mgw.laplace.integral} and \cite[(3.9)]{\BinghamDoney} with Cor.~\ref{c:moments.int}). Recalling \eqref{e:lapl.rem.expm},
$$\int_0^1
\f{\rem_{n,0}f(t;\vec v_t)}{t^\mmt}
	\f{dt}{t}
\le \f{1}{n!} \int_0^1
	\f{\E_\mgw[\abs{\vec Z_1}^n (1-e^{-t\abs{\vec Z_1}})]}
		{t^\be} \f{dt}{t}.$$
Applying Fubini's theorem and making the change of variable $t\mapsto \abs{\vec Z_1}t$ gives that the above is
$$= \f{1}{n!} \E_\mgw\Big[ \abs{\vec Z_1}^\mmt
	\int_0^{\abs{\vec Z_1}}
		\f{1-e^{-t}}{t^\be}
		\f{dt}{t}
	\Big]<\infty,$$
which concludes the proof.
\epf

\subsection{Harmonic moments and conductance estimates}
\label{ss:appx.cond}

In this section we prove the existence of harmonic moments for the normalized population size, extending part of \cite[Thm.~1]{\NeyVMoments} to the multi-type setting (using a similar proof). Using this result we adapt the methods of \cite[Lem.~2.2]{\PPPolar} to prove the conductance estimates used in the proofs of Cor.~\ref{c:timeroot} and Lem.~\ref{l:coupling.est}.

\blem\label{l:harmonic}
Assume \hypirred\ and \hypks. There exists some $r>0$ for which
$$\E_\mgw[W_\rt^{-r}\giv \nx^c]
\le \limsup_{n\to\infty} \E_\mgw[ \Znorm_n^{-r}\giv \nx^c]<\infty.$$

\bpf
Since the subtree of infinite descent described in \S\ref{ss:appx.gen} has the same normalized population size as the original tree, we may reduce to the case $\mgw(\nx)=0$. Expanding $F^{(n)}(\vec s)$ as a power series in $\vec s$ we find
$$
(F^{(n)}(\vec s))_a
\le \mgwtype{a}(\abs{\vec Z_n}=1) \, \nrm{\vec s}_\infty
+\mgwtype{a}(\abs{\vec Z_n}>1) \, \nrm{\vec s}_\infty^2.
$$
By \hypirred, there exists $n_0$ such that $\min_{a\in\cQ}\mgwtype{a}(\abs{\vec Z_n}>1)>0$ for all $n\ge n_0$, so that $F^{(n_0)}$ is a contraction on $[0,s_0]^\cQ$ for any $s_0<1$. By iterating this estimate, for any $s_0<1$ there exist constants $C<\infty$ and $\gam<1$ such that
$$\nrm{F^{(n)}(\vec s)}_\infty \le C\gam^n \nrm{\vec s}_\infty.$$
By Fubini's theorem,
$$\E_\mgw[\Znorm_n^{-r}]
= \f{\rho^{nr}}{\Gam(r)} \int_0^\infty f_n(u)\, u^{r-1} \d u,\quad
f_n(u)\equiv \E_\mgw[e^{-u \ip{\vec Z_n}{\rv}}].$$
We break up the integral into three parts, writing $I_{a,b}$ for the integral over $[a,b]$: By a change of variables,
$$\Gam(r)\,I_{0,\rho^{-n}}
= \int_0^1 \E_\mgw[e^{-u \Znorm_n}]
		\, u^{r-1} \d{u}
\le \f1r.$$
Next, we have
$$f_n(u)=\ip{\lv}{F^{(n)}(e^{-u\rvec{b}})}
\le \nrm{F^{(n)}(e^{-u\rvec{b}})}_\infty
\le C\gam^n \nrm{(e^{-u \rvec{b}})_{b\in\cQ}}_\infty
\le C\gam^n e^{-u \emin},$$
where for any $u_0>0$ we may choose constants $C<\infty$ and $\gam<1$ uniformly over all $u\ge u_0$. Therefore
$$\Gam(r)\,I_{1,\infty}
\le C(\rho^r\gam)^n
	\int_1^\infty e^{-u\emin} \, u^{r-1} \d u.$$
For $r>0$ small enough so that $\gam\rho^r<1$, we have $\lim_{n\to\infty} I_{1,\infty}=0$. It remains to consider
\begin{align*}
\Gam(r)\,I_{\rho^{-n},1}
&= \rho^{nr} \sum_{i=1}^n \int_{1/\rho^i}^{1/\rho^{i-1}}
       f_n(u) u^{r-1}\d{u}\\
&= \sum_{i=1}^n \rho^{r(n-i)} \int_1^\rho
       \E_\mgw[e^{-u\ip{\vec Z_n}{\rv}/\rho^i}]
        u^{r-1}\d{u}.
\end{align*}
By conditioning on the first $n-i$ levels of the tree,
\begin{align*}
&\E_\mgw[e^{-u\ip{\vec Z_n}{\rv}/\rho^i}]
= \E_\mgw\Big[
	\prod_{v\in D_{n-i}}
		\E_\mgw
		[e^{-u\ip{\vec Z^v_i}{\rv}/\rho^i} \giv \type_v]
	\Big] \\
&= \E_\mgw\Big[
	\prod_{a\in\cQ} \Phi^a_i(u)^{Z_{n-i}(a)}
	\Big]
=\ip{\lv}{ F^{(n-i)} [(\Phi^a_i(u))_{a\in\cQ}] },
\end{align*}
where $\Phi^a_i(u)\equiv\E_{\mgwtype{a}}[e^{-u\Znorm_i}]$.
But
$$\sup_{i\ge1} \sup_{u\ge1} \Phi^a_i(u)
= \sup_{i\ge1} \Phi^a_i(1)<1,$$
since $\Phi^a_i(1)<1$ for all $a,i$ and $\Phi^a_i(1) \to \E_{\mgwtype{a}}[e^{-u W_\rt}]$ which is less than $1$ by the Kesten--Stigum theorem as noted in \S\ref{ss:imgwr.wk.lim} (using \hypks). Therefore
$$\Gam(r)\,I_{\rho^{-n},1}
\le C\sum_{i=1}^n (\gam \rho^r)^{n-i} \int_1^\rho u^{r-1}\d u,$$
which is bounded in $n$ for small enough $r$. Putting the estimates together concludes the proof.
\epf
\elem

We conclude with the following conductance lower bound, a version of \cite[Lem.~2.2]{\PPPolar}. This clearly implies Lem.~\ref{l:polar} which was used in the proof of the annealed $\mgw$-CLT Propn.~\ref{p:mgw.aclt}.

\bppn\label{p:polar}
\bnm[(a)]
\item\label{p:polar.a}
Under \hypirred, \hypks, and \hypmmt{2}, there exist $0<r,C<\infty$ such that for all $\ep>0$,
$\mgw(\cond{\rt}{k}^{-1} \ge k^{1+\ep}\giv \nx^c) \le C k^{-r\ep}$.
\item\label{p:polar.q}
If further \hypmmt{\mmt} holds with $\mmt>2$, then for $\mgw$-a.e.\ $\tree\notin\nx$ there exists a random constant $C_\tree<\infty$ such that $\cond{\rt}{k}^{-1} \le C_\tree k$ for all $k$.
\enm

\bpf
(\ref{p:polar.a}) Recall that a \emph{unit flow} is a non-negative function $U$ on the vertices of $\tree$ such that for all $v\in\tree$, $U(v)=\sum_{w\in\desc v}U(w)$. For $v\in D_\ell$ define
$$U(v) =\f{W_v}{\sum_{u\in D_\ell} W_u} = \f{W_v}{\rho^\ell W_\rt};$$
it is easily seen that $U$ is a well-defined unit flow on $\nx^c$. It gives positive flow only to vertices of infinite descent, so by the discussion of \S\ref{ss:appx.gen} we may reduce to the case $\mgw(\nx)=0$. By Thomson's principle \cite[\S2.4]{lyons.peres}
$$\cond{\rt}{k}^{-1} \le \sum_{\ell=1}^k \rho^\ell \sum_{v\in D_\ell} U(v)^2
= \f{1}{W_\rt^2} \sum_{\ell=1}^k \f{1}{\rho^\ell} \sum_{v\in D_\ell} W_v^2.$$
By H\"older's inequality,
$$\E_\mgw[\cond{\rt}{k}^{-r}]
\le \E_\mgw \bigg[
	 \sum_{\ell=1}^k \f{1}{\rho^\ell} \sum_{v\in D_\ell} W_v^2
	\bigg]^r
	\E[W_\rt^{-2 r/(1-r)}]^{1-r} \le C k^r$$
for $r$ sufficiently small, using Lem.~\ref{l:harmonic} and $\mmt\ge2$. It follows from Markov's inequality that $\mgw(\cond{\rt}{k}^{-1} \ge k^{1+\ep}) \le C  k^{-r\ep}$.

\medskip\noindent
(\ref{p:polar.q}) We claim there exist $0<c,c'<\infty$ deterministic such that
\beq \label{e:w2.avg}
\mgw\bigg(
\f{1}{\abs{D_k}} \sum_{v\in D_k} W_v^2 \ge c' \bigg)
\le \rho^{-c k}
\eeq
Assuming the claim, we have
$$\rho^k \sum_{v\in D_k} U(v)^2
= \f{1}{\rho^k} \f{\sum_{v\in D_k} W_v^2}{W_\rt^2}
\le \f{C\Znorm_k}{W_\rt^2}
	\lp \f{1}{\abs{D_k}} \sum_{v\in D_k} W_v^2 \rp,$$
so by Borel--Cantelli
$$\limsup_{k\to\infty}\rho^k \sum_{v\in D_k} U(v)^2
\le \f{C c'}{W_\rt}<\infty,$$
which by Thomson's principle implies $\cond{\rt}{k}^{-1} \le C_\tree k$.

It remains to prove \eqref{e:w2.avg}. For any $1\le 1+r\le2 \wedge (\mmt/2)$, Lem.~\ref{l:petrov} and Markov's inequality give
\begin{align*}
&\mgw\bigg(
	\bigg|
	\f{1}{\abs{D_k}}
	\sum_{v\in D_k} (W_v^2-\E_{\mgwtype{\type_v}}[W_v^2])
	\bigg| \ge\ep
	\,\bigg|\, \cF_k
	\bigg)\\
&\le \f{2}{(\ep\abs{D_k})^{1+r}}
	\sum_{v\in D_k} \E_{\mgwtype{\type_v}}
		[\abs{W_v^2-\E_{\mgwtype{\type_v}}[W_v^2]}^{1+r}]
\le C(\ep,r) \abs{D_k}^{-r}.
\end{align*}
Taking expectations and applying Lem.~\ref{l:harmonic} then gives
$$\mgw\bigg(
	\bigg|
	\f{1}{\abs{D_k}}
	\sum_{v\in D_k} W_v^2
	-\f{1}{\abs{D_k}}
	\sum_{v\in D_k} \E_{\mgwtype{\type_v}}[W_v^2]
	\bigg| \ge\ep
	\bigg) \le  C(\ep,r) \rho^{-rk}$$
for $r$ sufficiently small. But
$$\f{1}{\abs{D_k}} \sum_{v\in D_k} \E_{\mgwtype{\type_v}}[W_v^2]$$
is clearly bounded uniformly in $k$ by a deterministic constant, so \eqref{e:w2.avg} is proved.
\epf
\eppn

\bibliography{mgw}
\bibliographystyle{abbrv}
\end{document}